\documentclass[10pt,a4paper]{article}

\usepackage[margin=1in]{geometry}
\usepackage{amsmath,amsthm,amssymb}
\allowdisplaybreaks
\usepackage[final]{pdfpages}
\usepackage{hyperref}
\usepackage{xfrac}
\usepackage{comment}
\usepackage{mathtools}

 \usepackage{array}
 \usepackage{ulem}
 \usepackage{verbatim}

\newcolumntype{M}[1]{>{\centering\arraybackslash}m{#1}}
\newcolumntype{N}{@{}m{0pt}@{}}

\newtheorem{theorem}{Theorem}[section]
\newtheorem{definition}[theorem]{Definition}
\newtheorem{lemma}[theorem]{Lemma}

\theoremstyle{remark}

\newtheorem{remark}[theorem]{Remark}

\newcommand{\N}{\mathbb{N}}
\newcommand{\Z}{\mathbb{Z}}
\newcommand{\R}{\mathbb{R}}

\title{Asymptotic preserving IMEX Active Flux schemes for a hyperbolic approximation of the Korteweg-de Vries equation}

\author{Simon Krotsch\thanks{Institute of Mathematics, University of W\"urzburg, Emil-Fischer-Stra\ss e 40, 97074 W\"urzburg, Germany, simon.krotsch@uni-wuerzburg.de (Simon Krotsch), christian.klingenberg@uni-wuerzburg.de (Christian Klingenberg)} , Wasilij Barsukow\thanks{Institut de Mathématiques de Bordeaux (IMB), CNRS UMR 5251, University of Bordeaux, 33405 Talence, France and Imperial College London, Abraham de Moivre Laboratory, CNRS IRL 2004, London, UK, wasilij.barsukow@math.u-bordeaux.fr} , Christian Klingenberg\footnotemark[1]}

\date{\vspace{-10mm}}

\begin{document}

\maketitle

\begin{abstract}
	We present a theoretical and numerical analysis of novel IMEX Active Flux schemes for a hyperbolic approximation of the Korteweg-de Vries (KdV) equation. For this approximation, it is important that the numerical schemes recover the limit equations (asymptotic preserving) and retain their order of accuracy in the limit (asymptotically accurate). Using formal asymptotic analysis and numerical experiments, we show that an Active Flux spatial discretisation (closely related to the standard semi-discrete Active Flux scheme), combined with an implicit-explicit Runge-Kutta scheme, is asymptotic preserving but not asymptotically accurate. We therefore propose alternative spatial discretisations that differ in the update equations of the point values. First, we derive a scheme by rescaling and eliminating terms in the upwind matrices. Second, we develop a second-order scheme by using central discretisation for the terms treated implicitly in time and present two modifications of it to third-order. Finally, we formulate a scheme with predefined upwind directions that are not chosen based on characteristics. We show that the alternative schemes are not only asymptotic preserving but also asymptotically accurate, and compare their accuracy, stability, and energy conservation numerically.
	
\end{abstract}

\section{Introduction}

In 1834, the Scottish engineer John Scott Russell observed a special type of wave while he was following a boat along the Union Channel: "[...] a rounded, smooth and well-defined heap of water, which continued its course along the channel apparently without change of form or diminution of speed."  \cite{russel_1845}. His observation and experiments to better comprehend the phenomena he encountered were followed by Boussinesq (1871) and Rayleigh (1876) and finally concluded in the derivation of a simple model partial differential equation in 1895 by Korteweg and de Vries \cite{korteweg_de_vries_1895} which admits solutions similar to Scott's observations. Today, this equation,
\begin{align}
	\frac{\partial}{\partial t}u(t,x) + u(t,x) \frac{\partial}{\partial x}u(t,x) + \frac{\partial^3}{\partial x^3}u(t,x) &= 0,
	\label{eq:kdv}
\end{align}
is known as the Korteweg-de Vries (KdV) equation. In fact, the KdV equation in a different form had already appeared in an article by Boussinesq in 1877, before it was rediscovered by Korteweg and de Vries \cite{marin_2012}, \cite{drazin_1983}. The special type of travelling wave solution is as follows:
\begin{theorem}[Single-soliton solution, \cite{evans_2010}]
	The two-parameter family of functions $u: \R^2 \rightarrow \R$  defined by
	\begin{equation}
		u(t,x) = \tilde u(x-ct) = 3c\ \mathrm{sech}^2 \left( \frac{\sqrt{c}}{2} (x-x_0-ct)\right) \label{eq:soliton}
	\end{equation}
	with $c > 0$, $x_0 \in \mathbb R$ satisfies the KdV equation for all $(t,x) \in \R^2$.
\end{theorem}
The parameter $c$ simultaneously controls the speed and amplitude of the solution. Due to the nonlinearity of \eqref{eq:kdv}, one cannot expect that the sum of solutions constitutes another solution. In 1965 Zabusky and Kruskal \cite{zabusky_kruskal} discovered that something close to a superposition principle holds true for travelling wave solutions of the kind of \eqref{eq:soliton}. Suppose that at $t=0$ two such solutions are given, separated from each other and the slower, smaller solution is located to the right of the larger, faster one. After sufficient time, the faster wave catches up to the slower and, up to a shift, their interaction looks as if the faster wave would simply overtake the slower. Zabusky and Kruskal referred to this type of special travelling wave solutions as \emph{solitons} \cite{ablowitz_segur_1981}. The following theorem states a solution to the KdV equation consisting of two solitons.
\begin{theorem}[Double-soliton solution, \cite{yoneyama_1984}]
	\label{th:doub_sol}
	The function $u:\R^2 \rightarrow \R$,
	\begin{equation*}
	u(t,x) = 12(\delta_2^2 - \delta_1^2) \frac{\delta_2^2 \mathrm{cosech}^2(\gamma_2(t,x)) +\delta_1^2 \mathrm{sech}^2(\gamma_1(t,x))}{(\delta_2\mathrm{coth}(\gamma_2(t,x)) -\delta_1\mathrm{tanh}(\gamma_1(t,x)))^2},
\end{equation*}
where $\gamma_i :\R^2 \rightarrow \R$, $\gamma_i(t,x) = \delta_i x - 4 \delta_i^3t$ for $i = 1,2$ and $\delta_1, \delta_2 \in \R$ with $\delta_1 < \delta_ 2$, satisfies the KdV equation for all $(t,x) \in \R^2$.
\end{theorem}
The KdV equation has an infinite number of conserved quantities (see \cite{biswas_etal_2025}). The first two are mass $M$ and energy $E$ with
\begin{align*}
	M(t) = \int_\R u(t,x) \mathrm{d}x &&\text{and}&& E(t) = \int_\R \frac{1}{2}(u(t,x))^2\mathrm{d}x.
\end{align*}

One can approximate high-order PDEs with a first-order hyperbolic system. This process is called hyperbolisation or hyperbolic relaxation. Hyperbolisation can serve multiple purposes, among them the easy implementation into existing codes and imposition of boundary conditions. In certain cases it can be used to remove numerical stiffness. Moreover, high-order PDEs allow for infinite propagation speed, which might be considered a violation of causality and is not the case for first order hyperbolic systems \cite{ketcheson_biswas_2025}. The hyperbolisations offer an approach to Riemann problems for PDEs that do not a priori admit discontinuous solutions. A substantial amount of research has been dedicated to hyperbolisation. For example, in recent years, hyperbolic approximations of the Navier-Stokes-Korteweg system \cite{dhaouadi_dumbser_2022}, the conduit equation \cite{gavrilyuk_etal_2024}, the Serre-Green-Naghdi model \cite{bassi_etal_2022}, and a class of dispersive and diffusive-dispersive equations \cite{barthwal_dhaouadi_rohde_2025} have been presented and analysed. The latter class includes the KdV equation.

Hyperbolisations of the same PDE can have different properties, for instance concerning energy conservation. In \cite{besse_etal_2022} the following hyperbolic approximation of the KdV equation was presented that was later also studied in \cite{biswas_etal_2025}:
\begin{subequations}
	\label{eq:hkdv}
	\begin{align}
	\frac{\partial}{\partial t}u(t,x) + u(t,x) \frac{\partial}{\partial x}u(t,x) + \frac{\partial}{\partial x}\psi(t,x) &= 0, \label{eq:hkdv_1} \\
	\frac{\partial}{\partial t}\rho(t,x) - \frac{1}{\varepsilon^2} \frac{\partial}{\partial x} \rho(t,x) &= - \frac{\psi(t,x)}{\varepsilon^2}, \label{eq:hkdv_2} \\
	\frac{\partial}{\partial t} \psi(t,x) + \frac{1}{\varepsilon^2}\frac{\partial}{\partial x} u(t,x) &= \frac{\rho(t,x)}{\varepsilon^2}, \label{eq:hkdv_3}
\end{align}
\end{subequations}
where $0 < \varepsilon <1 $ denotes the relaxation parameter. We use $\epsilon^2$ instead of $\epsilon$ to avoid its fractional powers in the matrices and the asymptotic analysis below. In the following, we refer to these equations as the hyperbolic KdV (HKdV) equations. They were originally introduced to ease the implementation of boundary conditions.
We can rewrite \eqref{eq:hkdv} as
\begin{equation}
	\frac{\partial}{\partial t}q(t,x) + \frac{\partial}{\partial x} f(q(t,x)) = Sq(t,x), \label{eq:balancelawgeneral}
\end{equation}
where
\begin{align*}
	q(t,x) = \begin{pmatrix}
		u(t,x) \\
		\rho(t,x) \\
		\psi (t,x)
	\end{pmatrix}, && f(q) := \begin{pmatrix}
		\frac{1}{2}(u)^2 + \psi  \\
		-\frac{1}{\varepsilon^2}\rho  \\
		\frac{1}{\varepsilon^2}u
	\end{pmatrix}, && \text{and} && S = \begin{pmatrix}
		0 & 0 & 0 \\
		0 & 0 & -\frac{1}{\varepsilon^2} \\
		0 & \frac{1}{\varepsilon^2} & 0
	\end{pmatrix}.
\end{align*}
The eigenvalues of the Jacobian of $f$ are
\begin{align}
	\label{eq:hkdv_eval}
	\lambda_1(q) = -\frac{1}{\varepsilon^2} , && \lambda_2(q) = \frac{1}{2}\left(u - \sqrt{u^2+\frac{4}{\varepsilon^2}} \right), && \lambda_3(q) = \frac{1}{2}\left(u + \sqrt{u^2+\frac{4}{\varepsilon^2}} \right)
\end{align}
and the eigenvectors are
\begin{align}
	\label{eq:hkdv_evec}
	v_1(q) = \begin{pmatrix}
		0 \\
		1 \\
		0
	\end{pmatrix}, && v_2(q) = \begin{pmatrix}
		\frac{\varepsilon^2}{2} \left( u - \sqrt{u^2 + \frac{4}{\varepsilon^2}}\right)\\
		0 \\
		1
	\end{pmatrix}, && v_3(q) = \begin{pmatrix}
		\frac{\varepsilon^2}{2}\left(u + \sqrt{u^2+\frac{4}{\varepsilon^2}} \right) \\
		0 \\
		1
	\end{pmatrix}.
\end{align}
Since the eigenvectors are linearly independent, the HKdV equations are a hyperbolic system of balance laws with a stiff linear source term.

 The formal limit of \eqref{eq:hkdv} is
\begin{subequations}
	\label{eq:lim_kdv}
	\begin{align}
	\frac{\partial}{\partial t}u(t,x) + u(t,x)\frac{\partial}{\partial x}u(t,x) + \frac{\partial^3}{\partial x^3}u(t,x) &= 0, \label{eq:lim_kdv_1}\\
	\rho(t,x) & = \frac{\partial}{\partial x}u(t,x) , \label{eq:lim_kdv_2} \\
	\psi(t,x) &= \frac{\partial^2}{\partial x^2}u(t,x). \label{eq:lim_kdv_3}
\end{align}
\end{subequations}
Equation \eqref{eq:lim_kdv_1} is the KdV equation \eqref{eq:kdv}, while $\rho$ and $\psi$ are the first and second derivatives of $u$. In \cite{besse_etal_2022} it was found that \eqref{eq:hkdv} also admits travelling wave solutions and in \cite{biswas_etal_2025} a way to compute them numerically is presented. In particular, there exist travelling wave solutions that tend to the single-soliton solution of the KdV equation (and its derivatives) as $\varepsilon \rightarrow 0$.

Just as the KdV equation, the HKdV equations possess conserved quantities.
For a finite $\varepsilon$, they conserve the mass $M$ and a modified energy
\begin{align*}
	E_\mathrm{H}(t) = \int_\R \frac{1}{2}\Big((u(t,x))^2 + \varepsilon^2(\rho(t,x))^2 + \varepsilon^2(\psi(t,x))^2 \Big) \mathrm{d}x.
\end{align*}
Upon $\varepsilon \rightarrow 0$, this energy converges formally to $E$.\\

\begin{figure}
	\centering
	\includegraphics[scale=0.5]{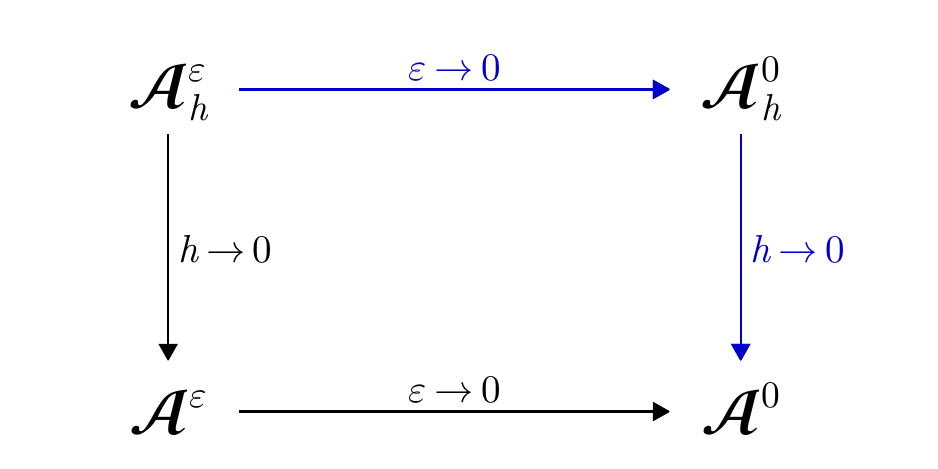}
	\caption{Illustration of the AP property.}
	\label{fig:ap_prop}
\end{figure}

Not every numerical method is able to obtain good approximations to the limit equations effectively by solving their hyperbolic approximations with a small $\epsilon$. Schemes that achieve this are called asymptotic preserving (AP).
This property is illustrated in Figure \ref{fig:ap_prop}. The hyperbolic approximation is $\mathcal{A}^\varepsilon$. We assume that this approximation features a relaxation parameter $\varepsilon$. If $\varepsilon \to 0$, $\mathcal{A}^\varepsilon$ formally becomes the limit equations denoted by $\mathcal{A}^0$ (under suitable assumptions on the initial data). $\mathcal{A}^\varepsilon_{h}$ is a (consistent and stable) numerical scheme discretising $\mathcal{A}^\varepsilon$, where $h$ denotes the discretisation parameter, e.g.\ $h = (\Delta x, \Delta t)$. For a fixed $h$, the asymptotic limit of $\mathcal{A}^\varepsilon_{h}$, provided it exists (under suitable assumptions on the discrete initial data), is denoted by $\mathcal{A}^0_{h}$.
\begin{definition}[Asymptotic preserving, \cite{jin_2012}, \cite{biswas_etal_2025}]
	\label{def:ap}
	If $\mathcal{A}^0_{h}$ is a consistent and stable discretisation of $\mathcal{A}^0$, then the numerical scheme $\mathcal{A}^\varepsilon_{h}$ is called asymptotic preserving (AP).
\end{definition}
If a numerical scheme is AP and, for example, able to achieve third order accuracy for sufficiently large $\varepsilon$, this does not immediately imply that the same holds as $\varepsilon \rightarrow 0$. Thus, it is natural to distinguish schemes that retain their order of accuracy as $\varepsilon \rightarrow 0$ from those that do not.
\begin{definition}[Asymptotically accurate, \cite{biswas_etal_2025}]
The numerical scheme $\mathcal{A}^\varepsilon_{h}$ is called asymptotically accurate (AA) if it retains its order of accuracy as $\varepsilon \rightarrow 0$.
\end{definition}
The limit equations of hyperbolic approximations often have the structure of \eqref{eq:lim_kdv} and consist of the (relevant) high-order PDE for $u$ and (irrelevant) constraints.
It is possible that a limit scheme $\mathcal{A}^0_{h}$ only discretises some part of the limit equations. Therefore, we call a scheme for the HKdV equations AP for the component $u$ if it includes a discretisation of \eqref{eq:lim_kdv_1}. The temporal and spatial discretisation scheme must be carefully chosen to ensure that the limit scheme $\mathcal{A}^0_{h}$ exists and is a discretisation of $\mathcal{A}^0$. Additionally, the temporal discretisation scheme, for a fixed $\Delta t$, has to remain stable as $\varepsilon \to 0$ (asymptotic stability).\\

In the early 2010s, Philip L. Roe introduced the Active Flux method \cite{eymann_roe_2011}, \cite{eymann_roe_2011_2}, \cite{eymann_roe_2013} as a new numerical method to solve hyperbolic conservation laws. It is based on a numerical method for the linear advection equation in one spatial dimension presented in 1977 by Bram van Leer \cite{van_leer_1977}. Roe was able to extend the idea of van Leer to multiple spatial dimensions. Today, Active Flux is a class of different methods. In addition to Roe's version of Active Flux, which is commonly referred to as the traditional or fully-discrete Active Flux method, there exists the generalised or semi-discrete Active Flux method \cite{abgrall_2023}, \cite{abgrall_barsukow_2023}, \cite{abgrall_barsukow_2023_2}, \cite{abgrall_barsukow_klingenberg_2025}. The different versions of Active Flux share some of the principles of the original scheme, notably (and contrary to classical finite volume or discontinuous Galerkin schemes), the Active Flux method computes the numerical solution without Riemann solvers. This is achieved by considering several types of degrees of freedom that are evolved in time independently: cell averages and point values at the cell interfaces in the case of the third order method. The point values are located at the interfaces of the cells and are shared between adjacent cells. This amounts to a globally continuous reconstruction.

The essential part of the fully-discrete Active Flux method is an evolution operator. It maps initial data to the exact or approximate solution of the problem at some given point in space and at some given time. For linear problems in one spatial dimension, it is straightforward to determine such an evolution operator, but finding good approximate evolution operators for sufficiently complicated problems in multiple spatial dimensions remains a significant challenge (with progress made in \cite{chudzik_helzel_2023}, \cite{barsukow_2026}, \cite{duraisamy_2026}, \cite{hensel_etal_2025}). In contrast, the semi-discrete Active Flux method does not depend on an evolution operator. Thus, a discretization determined once is immediately applicable to a variety of partial differential equations. For example, in \cite{duan_praveen_klingenberg_2025} the semi-discrete Active Flux method was used to numerically solve ideal magnetohydrodynamics.

The semi-discrete version itself offers a high degree of flexibility, with the particular choice of the degrees of freedom and the way to update them shaping the specific method. For instance, different ways to discretise the spatial derivative can be selected and, especially for systems, different ways to include upwinding can be considered. 

For the one-dimensional third-order semi-discrete Active Flux schemes considered here, we divide the real line into cells $\mathcal{C}_i = [x_{i-\frac{1}{2}}, x_{i+\frac{1}{2}}]$, where $x_{i \pm \frac{1}{2}} \in \R$ denotes the boundaries of the cell and $x_i$ the centre of the cell with $i \in \Z$. We assume that all cells have the same size $\Delta x = x_{i+\frac{1}{2}} - x_{i-\frac{1}{2}} > 0$. The degrees of freedom are cell averages and point values denoted by
\begin{align*}
	Q_i(t) \approx \frac{1}{\Delta x}\int_{\mathcal{C}_i}q(t, x) \mathrm{d}x &&\text{and}&& q_{i+\frac{1}{2}}(t) \approx q(t, x_{i+\frac{1}{2}}).
\end{align*}
Since the point values are shared between adjacent cells, this amounts to two degrees of freedom per cell: one point value and one average value. Suppose for a cell $\mathcal C_i$ the point values $q_{i \pm \frac{1}{2}}(t)$ and the cell average $Q_i(t)$ at time $t \geq 0$ are given. The quadratic polynomial
	\begin{align*}
		q_{\mathrm{recon},i} (t, \xi) = \frac{6Q_i(t) - q_{i+\frac{1}{2}}(t) - q_{i-\frac{1}{2}}(t)}{4} + \frac{q_{i+\frac{1}{2}}(t) - q_{i-\frac{1}{2}}(t)}{\Delta x}\xi +3 \frac{q_{i-\frac{1}{2}}(t) + q_{i+\frac{1}{2}}(t) - 2Q_i(t)}{\Delta x^2}\xi^2
	\end{align*}
	with $\xi \in [-\frac{\Delta x}{2}, \frac{\Delta x}{2}]$, interpolates the cell averages and the point values. Given this polynomial in each cell, we can define the globally continuous function $q_\mathrm{recon} : \R^+_0 \times \R \rightarrow \R^3$ by $q_\mathrm{recon}(t,x) = q_{\mathrm{recon},i} (t, x - x_i)$ if $x \in \mathcal C_i$. We call this function the (global) reconstruction. The time-continuous update equations of the point values and cell averages that constitute the specific schemes considered in this work are introduced in Sections \ref{sec:explicit_scheme}, \ref{sec:imex_scheme}, and \ref{sec:alternative_scheme}. 

It was found that the Active Flux method possesses structure preserving properties, e.g.\ it is stationarity preserving for the linear acoustic equations \cite{barsukow_etal_2026} and low Mach number compliant \cite{abgrall_barsukow_klingenberg_2025}. In \cite{duan_barsukow_klingenberg_2025_2}, the authors considered a hyperbolic approximation of the heat equation:
\begin{subequations}
	\label{eq:heat}
	\begin{align}
	\frac{\partial}{\partial t}p (t,\mathbf{x}) +  \frac{1}{\varepsilon} \nabla \cdot u(t,\mathbf{x}) &= 0, \\
	\frac{\partial}{\partial t} u(t,\mathbf{x}) + \frac{1}{\varepsilon} \nabla p(t,\mathbf{x})+\frac{\sigma}{\varepsilon^2}u(t,\mathbf{x}) &= 0,
\end{align}
\end{subequations}
where $\sigma >0$ and $0 < \varepsilon <1$. In the limit $\varepsilon \rightarrow 0$, this system formally becomes
\begin{align*}
	\frac{\partial}{\partial t} p(t,\mathbf{x})  &= \nabla \cdot \bigg( \frac{1}{\sigma} \nabla p(t,\mathbf{x})\bigg), \\
	u(t,\mathbf{x}) &= 0.
\end{align*}
They found that for \eqref{eq:heat}, the third-order semi-discrete Active Flux method with a point value update based on Jacobian splitting combined with implicit time discretisation is AP without any modification. In addition, in one spatial dimension they presented an alternative upwinding which resulted in an asymptotically accurate scheme. In this work we intend to further explore the ability of Active Flux to solve hyperbolic approximations of high-order partial differential equations. We are additionally interested in achieving not only the AP property, but to obtain a method that is AA.\\

To obtain asymptotically stable schemes, fully implicit time discretisation schemes can be used. However, in the case of nonlinear PDEs such as the HKdV equations, this amounts to solving nonlinear systems at each time step. A more economical approach is to discretise in time with implicit-explicit (IMEX) Runge-Kutta schemes. Consider an initial value problem for an ordinary differential equation,
\begin{align}
	\frac{\mathrm{d}}{\mathrm{d}t}v(t) &= f_1(v(t)) + f_2(v(t)) \label{eq:ode}
\end{align}
where $v:\R^+_0 \rightarrow \mathcal{O}$, $\mathcal{O} \subseteq \R^m$, $m\in \N$ and $f_1,f_2: \mathcal{O} \rightarrow \R^m$. Assume at time step $t^n \geq 0$, $n\in \N_0$, the approximation to $v(t^n)$, denoted by $v^n$, to be given. Denote by $v^{n+1}$ the approximations at the next time step $t^{n+1} = t^n + \Delta t$. An IMEX Runge-Kutta scheme with $\sigma \in \N$ stages for \eqref{eq:ode} takes the form
\begin{subequations}
\label{eq:2_4_1}
	\begin{align}
		v^{[k]} &= v^n + \Delta t \left(\sum_{j=1}^{k-1}\tilde a_{k,j}f_1(v^{[j]}) +  \sum_{j=1}^{\sigma}\hat a_{k,j}f_2(v^{[j]})\right), \label{eq:2_4_1_a}\\
		v^{n+1} &= v^n + \Delta t \left(\sum_{j=1}^{\sigma}\tilde b_{j}f_1(v^{[j]}) +  \sum_{j=1}^{\sigma}\hat b_{j}f_2(v^{[j]})\right),
	\end{align}
\end{subequations}
with $k=1,\dots \sigma$ and $v^{[k]}$ the solution at stage $k$. The lower triangular matrix $\tilde A = (\tilde a_{k,j}) \in \R^{\sigma \times \sigma}$ with zero diagonal and $\tilde b = (\tilde b_j) \in \R^\sigma$ are the coefficients of the explicit scheme, while the matrix $\hat A =  (\hat a_{k,j}) \in \R^{\sigma \times \sigma}$ and $\hat b = (\hat b_j) \in \R^\sigma$ are those of the implicit scheme. The matrices are such that the scheme is explicit in $f_1$ and implicit in $f_2$. We assume that $\hat A$ is a lower triangular matrix, i.e.\ the implicit part is a diagonally implicit Runge-Kutta scheme. This implies that the second sum in \eqref{eq:2_4_1_a} is from $1$ to $k$. 
\begin{definition}[Types of IMEX Runge-Kutta schemes, \cite{boscarino_russo_2024}]
	An IMEX Runge-Kutta scheme is of type I if the matrix $\hat A$ is invertible. It is said to be of type II if the matrix $\hat A$ can be written as
	\begin{align*}
		\hat A = \begin{pmatrix}
			0 & 0 \\
			a & A
		\end{pmatrix}
	\end{align*}
	with $a \in \R^{s}$, $s = \sigma -1$, and an invertible sub-matrix $A\in \R^{s\times s}$. If, in addition, $a = 0$ and $\hat b_1 = 0$, the IMEX Runge-Kutta scheme is called of type ARS.
\end{definition}
For IMEX Runge-Kutta schemes of type ARS, we can rewrite \eqref{eq:2_4_1} as
\begin{equation}
	\begin{aligned}
		v^{[0]} &= v^n, \\
		v^{[k]} &= v^n + \Delta t \left(\sum_{j=1}^{k}\tilde a_{k+1,j}f_1(v^{[j-1]}) +  \sum_{j=1}^{k} a_{k,j}f_2(v^{[j]})\right), \\
		v^{n+1} &= v^n + \Delta t \left(\sum_{j=1}^{\sigma}\tilde b_{j}f_1(v^{[j-1]}) +  \sum_{j=1}^{s} b_{j}f_2(v^{[j]})\right),
	\end{aligned}
	\label{eq:ars}
\end{equation}
where $k = 1, \dots s$ and $b = (\hat b_2, \dots, \hat b_\sigma)^\mathrm{T} \in \R^s$.
\begin{definition}[(Globally) stiffly accurate, \cite{ketcheson_biswas_2025}]
	An IMEX Runge-Kutta scheme is called globally stiffly accurate (GSA) if the corresponding diagonally implicit Runge-Kutta scheme is stiffly accurate, i.e.\ $\hat a_{\sigma j} = \hat b_j$ for $j = 1, \dots, \sigma$ and the explicit scheme is First Same As Last, i.e.\ $\tilde a_{\sigma,j} = \tilde b_j $ for all $j = 1, \dots \sigma$.
\end{definition}
An immediate consequence of this definition is the following: If an IMEX Runge-Kutta scheme is GSA, then $v^{n+1} = v^{[s]}.$ This property is useful for stiff problems \cite{ascher_etal_1997}.\\

In \cite{biswas_etal_2025}, the authors analysed IMEX Runge-Kutta time discretisations for the HKdV equations \eqref{eq:hkdv}. They found that an IMEX Runge-Kutta scheme of type I applied to a splitting of \eqref{eq:hkdv} without spatial discretisation is always AP for the $u$-component. If the scheme is additionally GSA, then the scheme is AP for all components. An IMEX Runge-Kutta scheme of type II applied to the same splitting is AP if the scheme is GSA and the initial data is well-prepared. In addition, they developed a structure preserving spatial discretisation based on summation-by-parts operators. If they applied IMEX Runge-Kutta schemes of type I and II to the discretisation instead of the space-continuous equations, they obtained the same results.

The main focus of this work is to analyse whether the Active Flux spatial discretisations are AP or AA. For this reason, we only consider IMEX Runge-Kutta schemes of type ARS. An additional reason for using these schemes is that their properties turn out to be important for the later analysis, e.g.\ the necessity for the initial data to be well-prepared can be avoided.
So far, \cite{duan_barsukow_klingenberg_2025_2} is the only study of AP properties of Active Flux for a hyperbolised high-order PDE. The CFL condition associated with the heat equation still allowed time-explicit methods to be used; the time-implicit method employed was fully implicit since the entire problem was linear. The present work, focusing on a nonlinear, third-order PDE must thus go significantly beyond \cite{duan_barsukow_klingenberg_2025_2}. Explicit methods can no longer be reasonably used, and the nonlinearity requires the development of an IMEX Active Flux method. This in itself is a novelty, since so far, the semi-discrete Active Flux scheme has only ever been integrated in time fully explicitly or fully implicitly.
The main contributions of this work are:
\begin{enumerate}
	\item[1.]  The development of an Active Flux scheme that can be discretised in time with IMEX Runge-Kutta schemes.
	\item[2.] The formal asymptotic analysis and numerical experiments showing that the novel IMEX Active Flux scheme is without any modification AP for the $u$ component and AP (for all variables) if the IMEX Runge-Kutta scheme is additionally GSA.
	\item[3.] The development and (theoretical and numerical) analysis of improved AP IMEX Active Flux schemes with alternative update equations of the point values, which, in particular, are also AA.
\end{enumerate}

The remainder of the paper is organised as follows: In Section \ref{sec:explicit_scheme}, the standard time-explicit Active Flux is recalled. In Section \ref{sec:imex_scheme}, a closely related scheme is derived by individually upwinding the nonlinear and stiff part of the HKdV equations. After that, a formal asymptotic analysis for this scheme discretised in time with IMEX Runge Kutta schemes of type $\mathrm{ARS}$ is conducted in Section \ref{sec:faa}. In Section \ref{sec:alternative_scheme}, alternative schemes are presented and analysed. Lastly, numerical results are presented in Section \ref{sec:num_exp} and the findings of this work are concluded in Section \ref{sec:con}.

\section{Active Flux scheme with coupled upwinding}
\label{sec:explicit_scheme}

We integrate the HKdV equations \eqref{eq:balancelawgeneral} over the cell $\mathcal C_i$ to obtain
\begin{align*}
	\frac{\mathrm{d}}{\mathrm{d}t} \int_{x_{i-\frac{1}{2}}}^{x_{i+\frac{1}{2}}} q(t,x) \mathrm{d}x + f(q(t,x_{i+\frac{1}{2}}))-f(q(t,x_{i-\frac{1}{2}}))= S\int_{x_{i-\frac{1}{2}}}^{x_{i+\frac{1}{2}}} q(t,x) \mathrm{d}x.
\end{align*}
The linear source term allows us to interchange the integral and matrix $S$. Based on this equation, the time-continuous update equation for the cell averages is
\begin{align}
	\frac{\mathrm d}{\mathrm{d}t} Q_i(t) = - \frac{f(q_{i+\frac{1}{2}}(t))) - f(q_{i-\frac{1}{2}}(t)))}{\Delta x} + SQ_i(t). \label{eq:av_exp}
\end{align}
To find an update equation for the point values, we consider the quasi-linear form \eqref{eq:hkdv} of the HKdV equations and approximate the spatial derivatives by finite difference formulas. We choose them to be high-order approximations of the derivatives at the location of the point values. Such formulas can be obtained by differentiating the reconstruction. At the location of a point value, there are two different derivatives: one from the cell on the right and one from the cell on the left.
We define
\begin{align*}
	D^-q_{i+\frac{1}{2}}(t) &= \frac{\partial}{\partial \xi} q_{\mathrm{recon},i+1} \bigg(t,-\frac{\Delta x}{2}\bigg) = \frac{1}{\Delta x}\left( -4q_{i+\frac{1}{2}}(t) + 6 Q_{i+1}(t) - 2q_{i+\frac{3}{2}}(t) \right)
\end{align*}
and
\begin{align*}
	D^+q_{i+\frac{1}{2}}(t) &= \frac{\partial}{\partial \xi} q_{\mathrm{recon},i} \bigg(t,\frac{\Delta x}{2}\bigg) = \frac{1}{\Delta x} \left( 2 q_{i-\frac{1}{2}}(t) - 6Q_i(t) + 4 q_{i+\frac{1}{2}}(t) \right).
\end{align*}
Having two different derivatives is not a disadvantage since it allows us to include upwinding. We diagonalise the Jacobian $F$ of $f$, $F(q) = R(q) \Lambda(q) R^{-1}(q)$ with $q \in \R^3$. $\Lambda(q)$ is a diagonal matrix with the eigenvalues \eqref{eq:hkdv_eval} on the diagonal. Next, we define $F^\pm(q) = R(q) \Lambda^\pm(q)R^{-1}(q)$ with $\Lambda^\pm (q) = \mathrm{diag}(\lambda_1^\pm(q),\lambda_2^\pm(q),\lambda_3^\pm(q))$, where $\lambda_i^+(q) = \max(\lambda_i(q),0)$ and $\lambda_i^-(q) = \min(\lambda_i(q),0)$.
The positive and negative parts of the Jacobian $F$ are:
\begin{align*}
	F^+(q) = \begin{pmatrix}
		\frac{\left(u+\sqrt{u^2+\frac{4}{\varepsilon^2}}\right)^2}{4 \sqrt{u^2+\frac{4}{\varepsilon^2}}} & 0 & \frac{1}{2} \left(\frac{u}{\sqrt{u^2+\frac{4}{\varepsilon^2}}}+1\right) \\
 0 & 0 & 0 \\
 \frac{1}{2\varepsilon^2} \left( \frac{u}{\sqrt{u^2+ \frac{4}{\varepsilon^2}}}+1\right) & 0 & \frac{1}{\varepsilon^2 \sqrt{u^2+\frac{4}{\varepsilon^2}}} \\
	\end{pmatrix}
\end{align*}
and
\begin{align*}
	F^-(q) = \begin{pmatrix}
		-\frac{\left(-u+\sqrt{u^2+\frac{4}{\varepsilon^2}}\right)^2}{4 \sqrt{u^2+\frac{4}{\varepsilon^2}}} & 0 & \frac{1}{2} \left(-\frac{u}{\sqrt{u^2+\frac{4}{\varepsilon^2}}}+1\right) \\
 0 & -\frac{1}{\varepsilon^2} & 0 \\
 \frac{1}{2\varepsilon^2} \left( -\frac{u}{\sqrt{u^2+ \frac{4}{\varepsilon^2}}}+1\right) & 0 & -\frac{1}{\varepsilon^2 \sqrt{u^2+\frac{4}{\varepsilon^2}}} \\
	\end{pmatrix}.
\end{align*}
The update equation of the point values is
\begin{align}
	\frac{\mathrm d}{\mathrm d t}q_{i+\frac{1}{2}}(t) = - \left(F^+(q_{i+\frac{1}{2}}(t)) D^+q_{i+\frac{1}{2}}(t) + F^-(q_{i+\frac{1}{2}}(t)) D^-q_{i+\frac{1}{2}}(t) \right) + Sq_{i+\frac{1}{2}}(t), \label{eq:pv_exp}
\end{align}
Observe that the Jacobian $F$ contains nonlinear and stiff linear terms. To distinguish the Active Flux scheme \eqref{eq:av_exp}, \eqref{eq:pv_exp} from other schemes developed in this work, it is referred to as the Active Flux scheme with combined upwinding or short the \texttt{combined} scheme. 

The \texttt{combined} scheme can be discretised  in time explicitly, for example, by a three-stage, third-order strong stability preserving Runge-Kutta (SSPRK3) scheme \cite{shu_osher_1988}. In this case, the time step size $\Delta t = t^{n+1} - t^n$ is subject to a CFL condition
\begin{align*}
	\Delta t = C_\mathrm{CFL}\frac{\Delta x}{\lambda_\mathrm{max}}, && \lambda_\mathrm{max} = \max_{i \in \{1,2,3\}} \max_{q \in \Gamma}|\lambda_i(q)|,
\end{align*}
where $\Gamma$ is the range of the solution $q$ at time $t^n$. In practice, $\Gamma$ is defined as the set of all point values at time step $t^n$. The eigenvalues depend on $\varepsilon$. As $\varepsilon \rightarrow 0$,  $|\lambda_i(q)| \rightarrow \infty$ for any $q \in \R^3$ and as a result $\Delta t \rightarrow 0$ in order for the scheme to remain stable, i.e.\ the scheme is not asymptotically stable and inefficient for small values of $\varepsilon$. Asymptotic stability can be obtained by discretising in time with IMEX Runge-Kutta schemes. Therefore, the update equations \eqref{eq:av_exp} and \eqref{eq:pv_exp} have to be split into a nonsingular and a linear singular part. The update equation of the cell averages \eqref{eq:av_exp} is in a form that makes splitting easily possible. It is more complicated to do the same for the update equation of the point values \eqref{eq:pv_exp} due to the presence of the square roots.

\section{Active Flux scheme with classical individual upwinding}
\label{sec:imex_scheme}
Inspired by \cite{biswas_etal_2025}, we follow a different approach to obtain an Active Flux scheme with update equations that are split into a nonsingular and a linear singular part. First, we split the spatial derivative in \eqref{eq:hkdv} into the two parts and then upwind them independently from each other. We write \eqref{eq:balancelawgeneral} as
\begin{align}
	\frac{\partial}{\partial t}q(t,x) + \frac{\partial}{\partial x} f_1(q(t,x)) + \frac{\partial}{\partial x} f_2(q(t,x)) = Sq(t,x),
	\label{eq:3_2_1}
\end{align}
where $q$ and $S$ are defined as usual and $f_1,f_2: \R^3 \rightarrow \R^3$ are
\begin{align*}
	f_1(q) = \begin{pmatrix}
		\frac{1}{2}u^2 \\[1pt]
		0\\[1pt]
		0
	\end{pmatrix} && \text{and} && f_2(q) = \begin{pmatrix}
		\psi \\[1pt]
		-\frac{1}{\varepsilon^2}\rho \\[1pt]
		\frac{1}{\varepsilon^2}u
	\end{pmatrix}
\end{align*}
The function $f_1$ is independent of $\varepsilon$ and $f_2$ is linear.
In quasi-linear form \eqref{eq:3_2_1} reads
\begin{align}
	\label{eq:qua_lin}
	\frac{\partial}{\partial t}q(t,x) + F_1(q(t,x))\frac{\partial}{\partial x} q(t,x) + F_2\frac{\partial}{\partial x} q(t,x) = Sq(t,x),
\end{align}
 where $F_1: \R^3 \rightarrow \R^{3 \times 3}$ with
\begin{align*}
	F_1(q)= \begin{pmatrix}
		u & 0 & 0 \\
		0 & 0 & 0 \\
		0 & 0 & 0
	\end{pmatrix} && \text{and} &&F_2 = \begin{pmatrix}
		0 & 0 & 1 \\
		0 & -\frac{1}{\varepsilon^2} & 0 \\
		\frac{1}{\varepsilon^2} & 0 & 0
	\end{pmatrix}
\end{align*}
are the Jacobians of $f_1$ and $f_2$. The Jacobian $F_2$ is independent of $q$. Hence, we dropped the dependence on it. Using the classical approach based on Jacobian splitting, we can determine the positive and negative parts of the two Jacobians. Starting with $F_1$, we have
\begin{align}
	F_1^+(q) = \begin{pmatrix}
		\max(u,0) & 0 & 0 \\
		0 & 0 & 0 \\
		0 & 0 & 0
	\end{pmatrix} &&\text{and} && F_1^-(q)  = \begin{pmatrix}
		\min(u,0) & 0 & 0 \\
		0 & 0 & 0 \\
		0 & 0 & 0
	\end{pmatrix}.
	\label{eq:jp_and_jm_1}
\end{align}
The positive and negative parts of $F_2$ are:
\begin{align}
	F_2^+ = \begin{pmatrix}
		\frac{1}{2\varepsilon} & 0 & \frac{1}{2} \\
		0 & 0 & 0 \\
		\frac{1}{2\varepsilon^2} & 0 & \frac{1}{2\varepsilon}
	\end{pmatrix} && \text{and} && F_2^- = \begin{pmatrix}
		-\frac{1}{2\varepsilon} & 0 & \frac{1}{2} \\
		0 & -\frac{1}{\varepsilon^2} & 0 \\
		\frac{1}{2\varepsilon^2} & 0 & -\frac{1}{2\varepsilon}
	\end{pmatrix}.
	\label{eq:jp_and_jm_2}
\end{align}
We define the update equations of the point values to be
\begin{align}
	 \begin{aligned}
		\frac{\mathrm{d}}{\mathrm{d}t}	q_{i+\frac{1}{2}}(t) = & - (F^+_1(q_{i+\frac{1}{2}}(t)) D^+ q_{i + \frac{1}{2}}(t) + F^-_1 (q_{i+\frac{1}{2}}(t))D^- q_{i + \frac{1}{2}}(t)) \\
         & - (F^+_2 D^+ q_{i + \frac{1}{2}}(t) + F^-_2 D^- q_{i + \frac{1}{2}}(t)) + Sq_{i + \frac{1}{2}}(t)
	\end{aligned}
	\label{eq:pv_cll}
\end{align}

These are not the same update equations as before since $F_1$ and $F_2$ are upwinded independently of each other. However, the resulting update equations are related: Expanding the square roots appearing in $F^\pm$ as power series in $\varepsilon$ one formally finds
\begin{align}
	F^+(q) = \begin{pmatrix}
		\frac{1}{2 \varepsilon} + \frac{u}{2} & 0 & \frac{1}{2} \\
		0 & 0 & 0 \\
		\frac{1}{2\varepsilon^2} + \frac{1}{4\varepsilon}u & 0 & \frac{1}{2\varepsilon}
	\end{pmatrix} + \mathcal{O}(\varepsilon) && \text{and} && F^-(q) = \begin{pmatrix}
		-\frac{1}{2 \varepsilon} + \frac{u}{2} & 0 & \frac{1}{2} \\
		0 & -\frac{1}{\varepsilon^2} & 0 \\
		\frac{1}{2\varepsilon^2} - \frac{1}{4\varepsilon}u & 0 & -\frac{1}{2\varepsilon}
	\end{pmatrix} + \mathcal{O}(\varepsilon).
	\label{eq:fp_and_fm}
\end{align}
Comparing \eqref{eq:fp_and_fm} to $F^\pm_1(q) +  F^\pm_2$, we observe that they differ by the additional terms $\pm \frac{1}{4\varepsilon} u$ in the third entry of the first row and by the $\mathcal{O}(\varepsilon)$ terms. 

From the integral form \eqref{eq:3_2_1}, we deduce in the same way as before the update equation of the cell averages:
\begin{align}
	\frac{\mathrm{d}}{\mathrm{d} t} Q_i(t) = -\frac{1}{\Delta x}\bigg(f_1(q_{i+\frac{1}{2}}(t)) - f_1(q_{i-\frac{1}{2}}(t)) \bigg) -\frac{1}{\Delta x}\bigg(f_2(q_{i+\frac{1}{2}}(t)) - f_2(q_{i-\frac{1}{2}}(t)) \bigg) + SQ_i(t)
	\label{eq:av_cll}
\end{align}
This is \eqref{eq:av_exp} split into a linear part and a nonlinear part, the latter being independent of the relaxation parameter.

Equations \eqref{eq:pv_cll} and \eqref{eq:av_cll} form what will be called the Active Flux scheme with classical individual upwinding or short the \texttt{classical} scheme. It allows us to discretise in time with IMEX Runge-Kutta schemes. However, it is still possible to discretise \eqref{eq:pv_cll} and \eqref{eq:av_cll} explicitly in time.

\section{Formal asymptotic analysis}
\label{sec:faa}
Next, we analyse whether the \texttt{classical} scheme discretised in time with an IMEX Runge-Kutta scheme \eqref{eq:ars} of type ARS is AP. To this end, we perform a formal asymptotic analysis, following the approaches of \cite{duan_barsukow_klingenberg_2025_2} and \cite{biswas_etal_2025}, and using Mathematica \cite{mathematica_2026} to assist in the analysis. We denote the point values and cell averages at time step $t^n$, $n\in \N_0$, by
\begin{align*}
	q_{i+\frac{1}{2}}^n = \begin{pmatrix}
		u_{i+\frac{1}{2}}^n \\
		\rho_{i+\frac{1}{2}}^n \\
		\psi_{i+\frac{1}{2}}^n
	\end{pmatrix} && \text{and} && Q_{i}^n = \begin{pmatrix}
		U_{i}^n \\
		P_{i}^n \\
		\Psi_{i}^n
	\end{pmatrix}
\end{align*}
with $i \in \Z$. The vector of stage solutions is defined as 
\begin{align}
	\overline \rho_{i+\frac{1}{2}} = \begin{pmatrix}
		\rho_{i+\frac{1}{2}}^{[1]} \\
		\vdots \\
		\rho_{i+\frac{1}{2}}^{[s]},
	\end{pmatrix} \label{eq:vectorstages}.
\end{align}
and $ \overline \psi_{i+\frac{1}{2}},\overline P_{i}$, and $\overline \Psi_{i}$ in a similar way. In addition, the temporal step size is denoted by $\Delta t = t^{n+1} - t^n$. The time discrete update equations for the cell averages are
\begin{subequations}
	\label{eq:tim_dis_av}
	\begin{align}
	U_{i}^{[k]} &=
	\begin{aligned}[t]
		&U^n_{i} + \Delta t \sum_{j=1}^{k}\tilde a_{k+1,j} \Biggl(-\frac{1}{\Delta x} \Biggl(\frac{1}{2}\left(u_{i+\frac{1}{2}}^{[j-1]}\right)^2 - \frac{1}{2}\left(u_{i-\frac{1}{2}}^{[j-1]}\right)^2 \Biggr)\Biggr) \\
		&+\Delta t \sum_{j=1}^k a_{k,j}\left(-\frac{1}{\Delta x} \left(\psi_{i+\frac{1}{2}}^{[j]} - \psi_{i-\frac{1}{2}}^{[j]}\right)\right),
	\end{aligned}\\
	U_{i}^{n+1} &= \begin{aligned}[t]
		&U^n_{i} + \Delta t \sum_{j=1}^{\sigma}\tilde b_{j} \Biggl(-\frac{1}{\Delta x} \Biggl(\frac{1}{2}\left(u_{i+\frac{1}{2}}^{[j-1]}\right)^2 - \frac{1}{2}\left(u_{i-\frac{1}{2}}^{[j-1]}\right)^2 \Biggr)\Biggr) \\
		&+\Delta t\sum_{j=1}^s b_{j}\left(-\frac{1}{\Delta x} \left(\psi_{i+\frac{1}{2}}^{[j]} - \psi_{i-\frac{1}{2}}^{[j]}\right)\right)
	\end{aligned}\\
	\overline P_{i} &=
	\begin{aligned}[t]
		&P^n_{i}\mathbb{I}_s + \Delta t A\left(-\frac{1}{\varepsilon^2}\frac{1}{\Delta x} \left(-\rho_{i+\frac{1}{2}}^{[j]} +\rho_{i-\frac{1}{2}}^{[j]}\right) - \frac{1}{\varepsilon^2} \Psi_i^{[j]}\right)_{j=1,\dots,s},
	\end{aligned} \label{eq:tim_dis_av_p}\\
	P_{i}^{n+1} &= \begin{aligned}[t]
		&P^n_{i}\mathbb{I}_s + \Delta t b \cdot \left(-\frac{1}{\varepsilon^2}\frac{1}{\Delta x} \left(-\rho_{i+\frac{1}{2}}^{[j]} +\rho_{i-\frac{1}{2}}^{[j]}\right) - \frac{1}{\varepsilon^2} \Psi_i^{[j]}\right)_{j=1,\dots,s},
	\end{aligned} \\
	\overline \Psi_{i} &=
	\begin{aligned}[t]
		&\Psi^n_{i}\mathbb{I}_s + \Delta t A \left(-\frac{1}{\varepsilon^2}\frac{1}{\Delta x} \left(u_{i+\frac{1}{2}}^{[j]} -u_{i-\frac{1}{2}}^{[j]}\right) + \frac{1}{\varepsilon^2} P_i^{[j]}\right)_{j=1,\dots,s},
	\end{aligned} \label{eq:tim_dis_av_psi}\\
	\Psi_{i}^{n+1} &= \begin{aligned}[t]
		&\Psi^n_{i}\mathbb{I}_s + \Delta t b \cdot \left(-\frac{1}{\varepsilon^2}\frac{1}{\Delta x} \left(u_{i+\frac{1}{2}}^{[j]} -u_{i-\frac{1}{2}}^{[j]}\right) + \frac{1}{\varepsilon^2} P_i^{[j]}\right)_{j=1,\dots,s}
	\end{aligned}
\end{align}
\end{subequations}
with $k=1,\dots s$ and $i \in \Z$. The ones for the point values are
\begin{subequations}
\label{eq:tim_dis_pv}
 	\begin{align}
	u_{i+\frac{1}{2}}^{[k]} &=
	\begin{aligned}[t]
		&u^n_{i+\frac{1}{2}} + \Delta t \sum_{j=1}^{k}\tilde a_{k+1,j} \left(-\max\left(u_{i+\frac{1}{2}}^{[j-1]},0\right)D^+u_{i+\frac{1}{2}}^{[j-1]}-\min\left(u_{i+\frac{1}{2}}^{[j-1]},0\right)D^-u_{i+\frac{1}{2}}^{[j-1]} \right) \\
		&+\Delta t \sum_{j=1}^k a_{k,j}\left(-\frac{1}{2 \varepsilon} \left(D^+u_{i+\frac{1}{2}}^{[j]} - D^-u_{i+\frac{1}{2}}^{[j]} \right) -  \frac{1}{2}\left(D^+\psi_{i+\frac{1}{2}}^{[j]} + D^-\psi_{i+\frac{1}{2}}^{[j]} \right) \right) ,
	\end{aligned}\label{eq:tim_dis_u}\\
	u_{i+\frac{1}{2}}^{n+1} &= \begin{aligned}[t]
		&u^n_{i+\frac{1}{2}} + \Delta t \sum_{j=1}^{\sigma}\tilde b_j \left(-\max\left(u_{i+\frac{1}{2}}^{[j-1]},0\right)D^+u_{i+\frac{1}{2}}^{[j-1]}-\min\left(u_{i+\frac{1}{2}}^{[j-1]},0\right)D^-u_{i+\frac{1}{2}}^{[j-1]} \right) \\
		&+\Delta t \sum_{j=1}^s b_{j}\left(-\frac{1}{2 \varepsilon} \left(D^+u_{i+\frac{1}{2}}^{[j]}- D^-u_{i+\frac{1}{2}}^{[j]} \right) -  \frac{1}{2}\left(D^+\psi_{i+\frac{1}{2}}^{[j]} + D^-\psi_{i+\frac{1}{2}}^{[j]} \right) \right)
	\end{aligned}\label{eq:tim_dis_u2} \\
	\overline \rho_{i+\frac{1}{2}} &=
	\begin{aligned}[t]
		&\rho^n_{i+\frac{1}{2}} \mathbb{I}_s + \Delta t A\left(\frac{1}{\varepsilon^2}D^-\rho_{i+\frac{1}{2}}^{[j]} - \frac{1}{\varepsilon^2} \psi_{i+\frac{1}{2}}^{[j]}\right)_{j=1,\dots,s},
	\end{aligned} \label{eq:tim_dis_pv_p}\\
	\rho_{i+\frac{1}{2}}^{n+1} &= \begin{aligned}[t]
		&\rho_{i+\frac{1}{2}}^n \mathbb{I}_s + \Delta t b \cdot \left(\frac{1}{\varepsilon^2}D^-\rho_{i+\frac{1}{2}}^{[j]} - \frac{1}{\varepsilon^2} \psi_{i+\frac{1}{2}}^{[j]}\right)_{j=1,\dots,s},
	\end{aligned} \label{eq:tim_dis_pv_p2} \\
	\overline \psi_{i+\frac{1}{2}} &=
	\begin{aligned}[t]
		&\psi^n_{i+\frac{1}{2}}\mathbb{I}_s + \Delta t A \bigg(-\frac{1}{2\varepsilon^2} \left(D^+u_{i+\frac{1}{2}}^{[j]} + D^-u_{i+\frac{1}{2}}^{[j]} \right) \\
        &- \frac{1}{2 \varepsilon} \left(D^+\psi_{i+\frac{1}{2}}^{[j]} - D^-\psi_{i+\frac{1}{2}}^{[j]} \right) + \frac{1}{\varepsilon^2} \rho_{i+\frac{1}{2}}^{[j]}\bigg)_{j=1,\dots,s} ,
	\end{aligned} \label{eq:tim_dis_pv_psi}\\
	\psi_{i+\frac{1}{2}}^{n+1} &= \begin{aligned}[t]
		&\psi_{i+\frac{1}{2}}^n\mathbb{I}_s + \Delta t b \cdot \bigg(-\frac{1}{2\varepsilon^2} \left(D^+u_{i+\frac{1}{2}}^{[j]} + D^-u_{i+\frac{1}{2}}^{[j]} \right) \\
        &- \frac{1}{2 \varepsilon} \left(D^+\psi_{i+\frac{1}{2}}^{[j]} - D^-\psi_{i+\frac{1}{2}}^{[j]} \right) + \frac{1}{\varepsilon^2} \rho_{i+\frac{1}{2}}^{[j]}\bigg)_{j=1,\dots,s},
	\end{aligned}\label{eq:tim_dis_pv_psi2}
\end{align}
\end{subequations}
The vector $\mathbb{I}_s \in \R^s$ in the above equations is defined as $\mathbb{I}_s  = (1,\dots,1)^\mathrm{T}$. To continue with the formal asymptotic analysis, we assume that the stage solutions for $j= 0, \dots, s$ and the numerical approximations at time $t^{n+1}$ are regular enough such that they can be expanded into power series in $\varepsilon$. For example,
\begin{align*}
	u_{i+\frac{1}{2}}^{[j]} &= u_{i+\frac{1}{2}}^{(0)[j]} + \varepsilon u_{i+\frac{1}{2}}^{(1)[j]} + \varepsilon^2 u_{i+\frac{1}{2}}^{(2)[j]} + \dots
\end{align*}
for $j = 0, \dots, s$. We now aim to determine the equations that the terms of order $0$ satisfy. Therefore we start by substituting the power series in $\varepsilon$ into the time-discrete update equation of the cell averages \eqref{eq:tim_dis_av} and ordering the resulting terms by the power of $\varepsilon$. This yields

\begin{subequations}
	\begin{align}
		0 &= \Delta t A\left(-\frac{1}{\Delta x} \left(-\rho_{i+\frac{1}{2}}^{(0)[j]} +\rho_{i-\frac{1}{2}}^{(0)[j]}\right) -  \Psi_i^{(0)[j]}\right)_{j=1,\dots,s}, \label{eq:av_faa_1_1}\\
		0 &= \Delta t b\cdot\left(-\frac{1}{\Delta x} \left(-\rho_{i+\frac{1}{2}}^{(0)[j]} +\rho_{i-\frac{1}{2}}^{(0)[j]}\right) -  \Psi_i^{(0)[j]}\right)_{j=1,\dots,s},
\end{align}
\end{subequations}
and
\begin{subequations}
	\begin{align}
		0 &= \Delta t A\left(-\frac{1}{\Delta x} \left(u_{i+\frac{1}{2}}^{(0)[j]} -u_{i-\frac{1}{2}}^{(0)[j]}\right) +  P_i^{(0)[j]}\right)_{j=1,\dots,s}, \label{eq:av_faa_2_1}\\
		0 &= \Delta t b\cdot\left(-\frac{1}{\Delta x} \left(u_{i+\frac{1}{2}}^{(0)[j]} -u_{i-\frac{1}{2}}^{(0)[j]}\right) +  P_i^{(0)[j]}\right)_{j=1,\dots,s},
\end{align}
\end{subequations}
and
\begin{subequations}
	\label{eq:3_3_3}
	 \begin{align}
		U_{i}^{(0)[k]} &=
	\begin{aligned}[t]
		&U^n_{i} + \Delta t \sum_{j=1}^{k}\tilde a_{k+1,j} \Biggl(-\frac{1}{\Delta x} \Biggl(\frac{1}{2}\left(u_{i+\frac{1}{2}}^{(0)[j-1]}\right)^2 - \frac{1}{2}\left(u_{i-\frac{1}{2}}^{(0)[j-1]}\right)^2 \Biggr)\Biggr) \\
		&+\Delta t \sum_{j=1}^k a_{k,j}\left(-\frac{1}{\Delta x} \left(\psi_{i+\frac{1}{2}}^{(0)[j]} - \psi_{i-\frac{1}{2}}^{(0)[j]}\right)\right),
	\end{aligned}\\
	U_{i}^{(0)n+1} &= \begin{aligned}[t]
		&U^n_{i} + \Delta t \sum_{j=1}^{\sigma}\tilde b_{j} \Biggl(-\frac{1}{\Delta x} \Biggl(\frac{1}{2}\left(u_{i+\frac{1}{2}}^{(0)[j-1]}\right)^2 - \frac{1}{2}\left(u_{i-\frac{1}{2}}^{(0)[j-1]}\right)^2 \Biggr)\Biggr) \\
		&+\Delta t\sum_{j=1}^s b_{j}\left(-\frac{1}{\Delta x} \left(\psi_{i+\frac{1}{2}}^{(0)[j]} - \psi_{i-\frac{1}{2}}^{(0)[j]}\right)\right)
	\end{aligned}
	\end{align}
\end{subequations}
We use the invertibility of $A$ to simplify \eqref{eq:av_faa_1_1} and \eqref{eq:av_faa_2_1} and obtain that the approximations of order $0$ satisfy
\begin{align}
	\Psi_i^{(0)[j]} &= \frac{1}{\Delta x}\left(\rho_{i+\frac{1}{2}}^{(0)[j]} -\rho_{i-\frac{1}{2}}^{(0)[j]}\right) \label{eq:3_3_5}
\end{align}
and
\begin{align}
	P_i^{(0)[j]} &= \frac{1}{\Delta x} \left(u_{i+\frac{1}{2}}^{(0)[j]} -u_{i-\frac{1}{2}}^{(0)[j]}\right)  \label{eq:3_3_6}
\end{align}
for all $j = 1, \dots s$. Equations \eqref{eq:3_3_5} and \eqref{eq:3_3_6} reassuringly look like discretisations of \eqref{eq:lim_kdv_2} and \eqref{eq:lim_kdv_3}. However, they cannot be used in \eqref{eq:3_3_3} to conclude that the latter is a discretisation of \eqref{eq:lim_kdv_1}, since \eqref{eq:3_3_3} contains the point values $\psi_{i+\frac12}$, and not the averages. We thus need to repeat the process for the update equations \eqref{eq:tim_dis_pv} of the point values before being able to proceed.

Inserting the power series into the time-discrete update equations of the point values \eqref{eq:tim_dis_pv} and ordering the resulting terms by the power of $\varepsilon$ gives
\begin{subequations}
	\begin{align}
		0 &= \Delta t A\left(D^-\rho_{i+\frac{1}{2}}^{(0)[j]} -  \psi_{i+\frac{1}{2}}^{(0)[j]}\right)_{j=1,\dots,s}, \label{eq:av_faa_3_1} \\
		0 &= \Delta t b \cdot \left(D^-\rho_{i+\frac{1}{2}}^{(0)[j]} -  \psi_{i+\frac{1}{2}}^{(0)[j]}\right)_{j=1,\dots,s}, \label{eq:av_faa_3_2}
\end{align}
\end{subequations}
and
\begin{subequations}
	\begin{align}
		0 &= \Delta t A\left(-\frac{1}{2} \left(D^+u_{i+\frac{1}{2}}^{(0)[j]} + D^-u_{i+\frac{1}{2}}^{(0)[j]} \right) +  \rho_{i+\frac{1}{2}}^{(0)[j]}\right)_{j=1,\dots,s},\label{eq:av_faa_4_1} \\
		0 &= \Delta t b \cdot \left(-\frac{1}{2} \left(D^+u_{i+\frac{1}{2}}^{(0)[j]} + D^-u_{i+\frac{1}{2}}^{(0)[j]} \right) +  \rho_{i+\frac{1}{2}}^{(0)[j]}\right)_{j=1,\dots,s}
\end{align}
\end{subequations}
and
\begin{subequations}
	\label{eq:3_3_1}
	\begin{align}
		u_{i+\frac{1}{2}}^{[k](0)} &=
	\begin{aligned}[t]
		&u^n_{i+\frac{1}{2}} + \Delta t \sum_{j=1}^{k}\tilde a_{k+1,j} \left(-\max\left((u_{i+\frac{1}{2}}^{(0)[j-1]},0\right)D^+u_{i+\frac{1}{2}}^{(0)[j-1]}\right.\\ 
		&\left. -\min\left(u_{i+\frac{1}{2}}^{(0)[j-1]},0\right)D^-u_{i+\frac{1}{2}}^{(0)[j-1]} \right) \\
		&+\Delta t \sum_{j=1}^k a_{k,j}\left(-\frac{1}{2} \left(D^+u_{i+\frac{1}{2}}^{(1)[j]} - D^-u_{i+\frac{1}{2}}^{(1)[j]} \right) -  \frac{1}{2}\left(D^+\psi_{i+\frac{1}{2}}^{(0)[j]} + D^-\psi_{i+\frac{1}{2}}^{(1)[j]} \right) \right),
	\end{aligned}\\
	u_{i+\frac{1}{2}}^{n+1(0)} &= \begin{aligned}[t]
		&u^n_{i+\frac{1}{2}} + \Delta t \sum_{j=1}^{\sigma}\tilde b_j \left(-\max\left(u_{i+\frac{1}{2}}^{[j-1](0)},0\right)D^+u_{i+\frac{1}{2}}^{[j-1](0)} \right.\\ 
		&\left. -\min\left(u_{i+\frac{1}{2}}^{[j-1](0)},0\right)D^-u_{i+\frac{1}{2}}^{[j-1](0)} \right) \\
		&+\Delta t \sum_{j=1}^s b_{j}\left(-\frac{1}{2} \left(D^+u_{i+\frac{1}{2}}^{(1)[j]}- D^-u_{i+\frac{1}{2}}^{(1)[j]} \right) -  \frac{1}{2}\left(D^+\psi_{i+\frac{1}{2}}^{(0)[j]} + D^-\psi_{i+\frac{1}{2}}^{(0)[j]} \right) \right)
	\end{aligned}
	\end{align}
\end{subequations}
We can again use the invertibility of $A$ to obtain
\begin{align}
	\psi_{i+\frac{1}{2}}^{(0)[j]}&= D^-\rho_{i+\frac{1}{2}}^{(0)[j]} \label{eq:3_3_7}
\end{align}
and
\begin{align}
	\rho_{i+\frac{1}{2}}^{(0)[j]} &= \frac{1}{2} \left(D^+u_{i+\frac{1}{2}}^{(0)[j]} + D^-u_{i+\frac{1}{2}}^{(0)[j]} \right) \label{eq:3_3_8}
\end{align}
for $i = 1\dots,s$ from \eqref{eq:av_faa_3_1} and \eqref{eq:av_faa_4_1}.
\begin{lemma}
	\label{lem:faa0}
	Assume $v_{i+\frac{1}{2}}$ and $V_{i}$ for $i \in \Z$ are point values and cell averages of an Active Flux scheme and $D^+v_{i+\frac{1}{2}}$ and $D^-v_{i+\frac{1}{2}}$ are defined as
	\begin{align*}
		D^+v_{i+\frac{1}{2}} &= \frac{1}{\Delta x} \left( 2 v_{i-\frac{1}{2}} - 6V_i + 4 v_{i+\frac{1}{2}} \right),\\
	 	D^-v_{i+\frac{1}{2}} &= \frac{1}{\Delta x}\left( -4v_{i+\frac{1}{2}} + 6 V_{i+1} - 2v_{i+\frac{3}{2}} \right).
	\end{align*}
	This implies
	\begin{subequations}
		\begin{align}
		\frac{1}{2}\left(D^+v_{i+\frac{1}{2}} + D^-v_{i+\frac{1}{2}}\right) &= \frac{1}{\Delta x}\left(v_{i-\frac{1}{2}} - 3 V_i + 3 V_{i+1} - v_{i+\frac{3}{2}}\right), \\
		\frac{1}{2}\left(D^+v_{i+\frac{1}{2}} - D^-v_{i+\frac{1}{2}}\right) &= \frac{1}{\Delta x}\left(v_{i-\frac{1}{2}} - 3 V_i +4v_{i+\frac{1}{2}}- 3 V_{i+1} + v_{i+\frac{3}{2}}\right). \label{eq:dp_m_dm}
	\end{align}
	\end{subequations}
	
\end{lemma}
\begin{proof}
	The result immediately follows from the definition of $D^+v_{i+\frac{1}{2}}$ and $D^-v_{i+\frac{1}{2}}$.
\end{proof}
\begin{lemma}
	\label{lem:faa1}
	The terms of order $0$ satisfy
	\begin{equation*}
		\begin{aligned}
			\psi_{i+\frac{1}{2}}^{(0)[j]} =  \frac{2}{\Delta x^2} \left(-2u_{i-\frac{1}{2}}^{(0)[j]} + 6U_i^{(0)[j]}- 4u_{i+\frac{1}{2}}^{(0)[j]} - 3 U_{i+1}^{(0)[j]} +5 u_{i+\frac{3}{2}}^{(0)[j]} - 3 U_{i+2}^{(0)[j]} + u_{i+\frac{5}{2}}^{(0)[j]}\right)
		\end{aligned}
	\end{equation*}
	and
	\begin{equation*}
		\frac{1}{2}\left(D^+\psi_{i+\frac{1}{2}}^{(0)[j]} + D^-\psi_{i+\frac{1}{2}}^{(0)[j]} \right) = \begin{aligned}[t]
			& \frac{1}{\Delta x^3}\bigg(-u_{i-\frac{3}{2}}^{(0)[j]} + 3U_{i-1}^{(0)[j]} -14u_{i-\frac{1}{2}}^{(0)[j]} + 21U_i^{(0)[j]} + 14u_{i+\frac{1}{2}}^{(0)[j]} \\
			& -45U_{i+1}^{(0)[j]}  +16 u_{i+\frac{3}{2}}^{(0)[j]} +15U_{i+2}^{(0)[j]} -13u_{i+\frac{5}{2}}^{(0)[j]} + 6U_{i+3}^{(0)[j]} -2 u_{i+\frac{7}{2}}^{(0)[j]}\bigg)
		\end{aligned}
	\end{equation*}
	for $j = 1, \dots, s$.
\end{lemma}
\begin{proof}
Using the Equations \eqref{eq:3_3_5}, \eqref{eq:3_3_6}, \eqref{eq:3_3_7}, and \eqref{eq:3_3_8} as well as Lemma \ref{lem:faa0}, we obtain
\begin{align*}
	\psi_{i+\frac{1}{2}}^{(0)[j]} &\overset{\eqref{eq:3_3_7}}{=} D^-\rho_{i+\frac{1}{2}}^{(0)[j]} \\
	&= \frac{1}{\Delta x}\left(-4\rho_{i+\frac{1}{2}}^{(0)[j]} + 6P_{i+1}^{(0)[j]} -2 \rho_{i+\frac{3}{2}}^{(0)[j]}\right) \\
	&\!\!\!\!\!\!\overset{\eqref{eq:3_3_8},\eqref{eq:3_3_6}}{=} \begin{aligned}[t]
		&\frac{1}{\Delta x} \left( - 2\left(D^+u_{i+\frac{1}{2}}^{(0)[j]} + D^-u_{i+\frac{1}{2}}^{(0)[j]}\right) +\frac{6}{\Delta x} \left(u_{i+\frac{3}{2}}^{(0)[j]} - u_{i+\frac{1}{2}}^{(0)[j]}\right) \right.\left. - \left(D^+u_{i+\frac{3}{2}}^{(0)[j]} + D^-u_{i+\frac{3}{2}}^{(0)[j]}\right)\right)
	\end{aligned}\\
	&= \frac{2}{\Delta x^2} \left(-2u_{i-\frac{1}{2}}^{(0)[j]} + 6U_i^{(0)[j]}- 4u_{i+\frac{1}{2}}^{(0)[j]} - 3 U_{i+1}^{(0)[j]} +5 u_{i+\frac{3}{2}}^{(0)[j]} - 3 U_{i+2}^{(0)[j]} + u_{i+\frac{5}{2}}^{(0)[j]}\right).
\end{align*}
In a similar way, using Lemma \ref{lem:faa0}, we can deduce
\begin{align*}
	\frac{1}{2}&\left(D^+\psi_{i+\frac{1}{2}}^{(0)[j]} + D^-\psi_{i+\frac{1}{2}}^{(0)[j]}\right) = \ \frac{1}{\Delta x}\left(\psi_{i-\frac{1}{2}}^{(0)[j]} - 3 \Psi_i^{(0)[j]} + 3 \Psi_{i+1}^{(0)[j]} - \psi_{i+\frac{3}{2}}^{(0)[j]}\right)\\
	&\!\!\!\!\!\!\overset{\eqref{eq:3_3_7},\eqref{eq:3_3_5}}{=} \begin{aligned}[t]
		&\frac{1}{\Delta x}\left(D^-\rho_{i-\frac{1}{2}}^{(0)[j]} -\frac{3}{\Delta x} \left(\rho_{i+\frac{1}{2}}^{(0)[j]} - \rho_{i-\frac{1}{2}}^{(0)[j]}\right) + \frac{3}{\Delta x} \left(\rho_{i+\frac{3}{2}}^{(0)[j]} - \rho_{i+\frac{1}{2}}^{(0)[j]}\right) - D^-\rho_{i+\frac{3}{2}}^{(0)[j]} \right)
	\end{aligned}\\
	&= \frac{1}{\Delta x^2} \left(-\rho_{i-\frac{1}{2}}^{(0)[j]} + 6P_i^{(0)[j]} -8 \rho_{i+\frac{1}{2}}^{(0)[j]} +7  \rho_{i+\frac{3}{2}}^{(0)[j]} -6P_{i+2}^{(0)[j]}  +2\rho_{i+\frac{5}{2}}^{(0)[j]}\right) \\
	&\!\!\!\!\!\!\overset{\eqref{eq:3_3_8},\eqref{eq:3_3_6}}{=}\begin{aligned}[t]
		 &\frac{1}{\Delta x^2} \left(-\frac{1}{2}\left(D^+u_{i-\frac{1}{2}}^{(0)[j]} + D^-u_{i-\frac{1}{2}}^{(0)[j]}\right) +\frac{6}{\Delta x} \left(u_{i+\frac{1}{2}}^{(0)[j]} - u_{i-\frac{1}{2}}^{(0)[j]}\right) -4\left(D^+u_{i+\frac{1}{2}}^{(0)[j]} + D^-u_{i+\frac{1}{2}}^{(0)[j]}\right) 
		 \right.\\& \left.
		 + \frac{7}{2}\left(D^+u_{i+\frac{3}{2}}^{(0)[j]} + D^-u_{i+\frac{3}{2}}^{(0)[j]}\right) -\frac{6}{\Delta x} \left(u_{i+\frac{5}{2}}^{(0)[j]} - u_{i+\frac{3}{2}}^{(0)[j]}\right) +\left(D^+u_{i+\frac{5}{2}}^{(0)[j]} + D^-u_{i+\frac{5}{2}}^{(0)[j]}\right)\right)
	\end{aligned}\\
	&= \begin{aligned}[t]
		&\frac{1}{\Delta x^3}\left(-u_{i-\frac{3}{2}}^{(0)[j]} + 3U_{i-1}^{(0)[j]} -14u_{i-\frac{1}{2}}^{(0)[j]} + 21U_i^{(0)[j]} + 14u_{i+\frac{1}{2}}^{(0)[j]}  \right.\\
	& \left.-45U_{i+1}^{(0)[j]} +16 u_{i+\frac{3}{2}}^{(0)[j]} +15U_{i+2}^{(0)[j]} -13u_{i+\frac{5}{2}}^{(0)[j]} + 6U_{i+3}^{(0)[j]} -2 u_{i+\frac{7}{2}}^{(0)[j]}\right).
	\end{aligned}
\end{align*}
\end{proof}
Lemma \ref{lem:faa1} implies that \eqref{eq:3_3_1} can be understood as
\begin{equation}
	\begin{aligned}
		\frac{\mathrm{d}}{\mathrm{d}t} u_{i+\frac{1}{2}}^{(0)}(t) = \begin{aligned}[t]
			& -\bigg(\max\left(u_{i+\frac{1}{2}}^{(0)}(t) ,0\right) \frac{1}{\Delta x}\left(2u_{i-\frac{1}{2}}^{(0)}(t) - 6U_i^{(0)}(t) + 4 u_{i+\frac{1}{2}}^{(0)}(t)\right) \\
			& + \min\left(u_{i+\frac{1}{2}}^{(0)}(t) ,0\right) \frac{1}{\Delta x}\left(-4u_{i+\frac{1}{2}}^{(0)}(t) + 6 U_{i+1}^{(0)}(t) - 2 u_{i+\frac{3}{2}}^{(0)}(t)  \right) \bigg)  \\
			& -\frac{1}{\Delta x^3}\bigg(-u_{i-\frac{3}{2}}^{(0)}(t) + 3U_{i-1}^{(0)}(t) -14u_{i-\frac{1}{2}}^{(0)}(t) + 21U_i^{(0)}(t) + 14u_{i+\frac{1}{2}}^{(0)}(t) \\
			& -45U_{i+1}^{(0)}(t)  +16 u_{i+\frac{3}{2}}^{(0)}(t) +15U_{i+2}^{(0)}(t) -13u_{i+\frac{5}{2}}^{(0)}(t) + 6U_{i+3}^{(0)}(t) -2 u_{i+\frac{7}{2}}^{(0)}(t)\bigg) \\
			& -\frac{1}{\Delta x}\left(u_{i-\frac{1}{2}}^{(1)}(t) - 3 U_i^{(1)}(t) +4u_{i+\frac{1}{2}}^{(1)}(t)- 3 U_{i+1}^{(1)}(t) + u_{i+\frac{3}{2}}^{(1)}(t)\right)
		\end{aligned}
	\end{aligned}
	\label{eq:3_3_2}
\end{equation}
discretised using an IMEX Runge-Kutta scheme and that, similarly, \eqref{eq:3_3_3} is the temporal discretisation of
\begin{equation}
	\begin{aligned}
	\frac{\mathrm{d}}{\mathrm{d}t}U_i^{(0)}(t) = \begin{aligned}[t]
			& -\frac{1}{\Delta x}\left(\frac{1}{2}\left(u_{i+\frac{1}{2}}^{(0)}(t)\right)^2 -\frac{1}{2}\left(u_{i-\frac{1}{2}}^{(0)}(t)\right)^2  \right)  \\
			&- \frac{2}{\Delta x^3}\left(\bigg(-2u_{i-\frac{1}{2}}^{(0)}(t) + 6U_i^{(0)}(t)- 4u_{i+\frac{1}{2}}^{(0)}(t) - 3 U_{i+1}^{(0)}(t) +5 u_{i+\frac{3}{2}}^{(0)}(t) \right. \\
			&\left.- 3 U_{i+2}^{(0)}(t) + u_{i+\frac{5}{2}}^{(0)}(t)\right) - \left(-2u_{i-\frac{3}{2}}^{(0)}(t) + 6U_{i-1}^{(0)}(t)- 4u_{i-\frac{1}{2}}^{(0)}(t) - 3 U_{i}^{(0)} (t) \right.\\
			&\left.\left.+5 u_{i+\frac{1}{2}}^{(0)}(t) - 3 U_{i+1}^{(0)}(t) + u_{i+\frac{3}{2}}^{(0)}(t)\right)\right)
		\end{aligned}
	\end{aligned}
	\label{eq:3_3_4}
\end{equation}
with an error in $\mathcal{O}(\Delta t^p)$, where $p$ denotes the order of the IMEX Runge-Kutta scheme. Observe how we were able to recover time-continuous update equations for $u$ (point values and cell averages) that do not contain point values or cell averages of $\rho$ and $\psi$. If the spatial discretisations appearing therein are approximations of the spatial derivatives in \eqref{eq:lim_kdv_1}, then the \texttt{classical} scheme (discretised in time with an IMEX-Runge Kutta scheme of type ARS) is AP for the component $u$. For the scheme to be AP for all components, we need to recover discretisations of \eqref{eq:lim_kdv_2} and \eqref{eq:lim_kdv_3}. Equations \eqref{eq:3_3_5}, \eqref{eq:3_3_6}, \eqref{eq:3_3_7}, and \eqref{eq:3_3_8} look like discretisations of \eqref{eq:lim_kdv_2} and \eqref{eq:lim_kdv_3} but they only hold for stage solutions and not for the approximation at the time steps. If at this point, we assume that the IMEX Runge-Kutta scheme is GSA, then the last stage solution is equal to the numerical approximation at the time step $t^{n+1}$. Using Lemma \ref{lem:faa1} once more, we then obtain 
\begin{align}
	\rho_{i+\frac{1}{2}}^{(0)n+1} &=  \frac{1}{\Delta x} \left(u_{i-\frac{1}{2}}^{(0)n+1} - 3U_i^{(0)n+1} + 3U_{i+1}^{(0)n+1} - u_{i+\frac{3}{2}}^{(0)n+1}  \right) \label{eq:3_3_11}
\end{align}
and
\begin{equation}
		\psi_{i+\frac{1}{2}}^{(0)n+1} = \begin{aligned}[t]
			&\frac{2}{\Delta x^2} \left(-2u_{i-\frac{1}{2}}^{(0)n+1} + 6U_i^{(0)n+1} - 4u_{i+\frac{1}{2}}^{(0)n+1}\right.\\
		&\left. - 3 U_{i+1}^{(0)n+1} +5 u_{i+\frac{3}{2}}^{(0)n+1} - 3 U_{i+2}^{(0)n+1} + u_{i+\frac{5}{2}}^{(0)n+1}\right).
		\end{aligned} 
	\label{eq:3_3_12}
\end{equation}
and
\begin{align}
	P_{i}^{(0)n+1}  = \frac{1}{\Delta x} \left(u_{i+\frac{1}{2}}^{(0)n+1} - u_{i-\frac{1}{2}}^{(0)n+1} \right)
	  \label{eq:3_3_9}
\end{align}
and
\begin{equation}
	\begin{aligned}
		\Psi_i^{(0)n+1} =\begin{aligned}[t] &\frac{1}{\Delta x^2} \left(\left(u_{i-\frac{1}{2}}^{(0)n+1} - 3U_i^{(0)n+1} + 3U_{i+1}^{(0)n+1} - u_{i+\frac{3}{2}}^{(0)n+1}  \right)\right. \\
		&- \left.\left(u_{i-\frac{3}{2}}^{(0)n+1} - 3U_{i-1}^{(0)n+1} + 3U_{i}^{(0)n+1} - u_{i+\frac{1}{2}}^{(0)n+1}  \right)  \right). \end{aligned}
	\end{aligned}\label{eq:3_3_10}
\end{equation}
These four equations are temporal discretisations of the following time-continuous equations:
\begin{align}
	\rho_{i+\frac{1}{2}}^{(0)}(t) &= \frac{1}{\Delta x} \left(u_{i-\frac{1}{2}}^{(0)}(t) - 3U_i^{(0)}(t) + 3U_{i+1}^{(0)}(t) - u_{i+\frac{3}{2}}^{(0)}(t)  \right)
	\label{eq:3_3_15}
	\end{align}
and
\begin{equation}
	\begin{aligned}
	\psi_{i+\frac{1}{2}}^{(0)}(t)  = \begin{aligned}[t]
		 &\frac{2}{\Delta x^2} \left(-2u_{i-\frac{1}{2}}^{(0)}(t)  + 6U_i^{(0)}(t)  - 4u_{i+\frac{1}{2}}^{(0)}(t)  - 3 U_{i+1}^{(0)}(t) \right. \left. +5 u_{i+\frac{3}{2}}^{(0)}(t)  - 3 U_{i+2}^{(0)}(t)  + u_{i+\frac{5}{2}}^{(0)}(t) \right)
		 \label{eq:3_3_16}
	\end{aligned}
\end{aligned}
\end{equation}
and
\begin{align}
	P_{i}^{(0)}(t) &= \frac{1}{\Delta x} \left(u_{i+\frac{1}{2}}^{(0)}(t) - u_{i-\frac{1}{2}}^{(0)}(t) \right)  \label{eq:3_3_13}
\end{align}
and
\begin{equation}
	\begin{aligned}
		\Psi_i^{(0)}(t) = \begin{aligned}[t]
		&\frac{1}{\Delta x^2} \left(\left(u_{i-\frac{1}{2}}^{(0)}(t) - 3U_i^{(0)}(t) + 3U_{i+1}^{(0)}(t) - u_{i+\frac{3}{2}}^{(0)}(t)  \right) \right. \\
		 &\left.- \left(u_{i-\frac{3}{2}}^{(0)}(t) - 3U_{i-1}^{(0)}(t) + 3U_{i}^{(0)}(t) - u_{i+\frac{1}{2}}^{(0)}(t)  \right)  \right).
		\end{aligned}
	\end{aligned}
	\label{eq:3_3_14}
\end{equation}
If \eqref{eq:3_3_15} and \eqref{eq:3_3_13} are (spatial) discretisations of \eqref{eq:lim_kdv_2}, and if \eqref{eq:3_3_16} and \eqref{eq:3_3_14} are discretisations of \eqref{eq:lim_kdv_3}, then the scheme is AP. This is indeed the case and the next theorem establishes the AP property.
\begin{theorem}
	\label{th:3_3_1}
	Assume $u^{(0)}, u^{(1)}: \R_0^+ \times \R \rightarrow \R$ are regular functions with bounded derivatives that interpolate the point values and cell averages. Under this assumption, \eqref{eq:3_3_2} and \eqref{eq:3_3_4} are, respectively, discretisations of 
	\begin{align}
		\label{eq:th_1}
		\frac{\partial}{\partial t}u^{(0)}(t,x) + u^{(0)}(t,x) \frac{\partial}{\partial x}u^{(0)}(t,x) + \frac{\partial^3}{\partial x^3}u^{(0)}(t,x) = 0
	\end{align}
	and
	\begin{align}
		\label{eq:th_2}
		\frac{1}{\Delta x}\int_{\mathcal{C}_i} \left( \frac{\partial}{\partial t}u^{(0)}(t,x) + u^{(0)}(t,x) \frac{\partial}{\partial x}u^{(0)}(t,x) + \frac{\partial^3}{\partial x^3} u^{(0)}(t,x) \right ) \mathrm{d}x = 0
	\end{align}
	with an error $\mathcal{O}(\Delta x)$. In addition, \eqref{eq:3_3_15} and \eqref{eq:3_3_13} are, respectively, discretisations of 
	\begin{align}
		\label{eq:th_3}
		\rho^{(0)}(t,x) &= \frac{\partial}{\partial x}u^{(0)}(t,x)
	\end{align}
	 with an error $\mathcal{O}(\Delta x^2)$ and
	\begin{align}
		\label{eq:th_4}
		\frac{1}{\Delta x}\int_{\mathcal{C}_i} \rho^{(0)}(t,x)\mathrm{d}x &= \frac{1}{\Delta x}\int_{\mathcal{C}_i} \frac{\partial}{\partial x}u^{(0)}(t,x) \mathrm{d}x.
	\end{align}
	with no discretization error. Equation \eqref{eq:3_3_16} is a discretisation of
	 \begin{align}
	 	\label{eq:th_5}
		\psi^{(0)}(t,x) &= \frac{\partial^2}{\partial x^2}u^{(0)}(t,x)
	\end{align}
	  with an error $\mathcal{O}(\Delta x)$, while \eqref{eq:3_3_14} discretises
\begin{align}
	\label{eq:th_6}
	\frac{1}{\Delta x}\int_{\mathcal{C}_i} \psi^{(0)}(t^,x)\mathrm{d}x &= \frac{1}{\Delta x}\int_{\mathcal{C}_i} \frac{\partial^2}{\partial x^2}u^{(0)}(t,x) \mathrm{d}x
\end{align}
with an error $\mathcal{O}(\Delta x^2)$.
\end{theorem}
\begin{proof} 
Starting with \eqref{eq:3_3_2}, it follows that 
\begin{multline*}
	\bigg(\max\left(u_{i+\frac{1}{2}}^{(0)}(t),0\right)D^+u_{i+\frac{1}{2}}^{(0)}(t) + \min\left(u_{i+\frac{1}{2}}^{(0)}(t),0\right) D^-u_{i+\frac{1}{2}}^{(0)}(t))\bigg) \\
	= u^{(0)} (t,x_{i+\frac{1}{2}})\frac{\partial}{\partial x}u^{(0)}(t,x_{i+\frac{1}{2}}) - u^{(0)} (t,x_{i+\frac{1}{2}})\frac{\Delta x^2}{12}\frac{\partial^3}{\partial x^3}u^{(0)}(t,x_{i+\frac{1}{2}}) + \mathcal{O}(\Delta x^3)
\end{multline*}
and 
\begin{multline*}
		\frac{1}{\Delta x^3}\left(-u_{i-\frac{3}{2}}^{(0)}(t) + 3U_{i-1}^{(0)}(t) -14u_{i-\frac{1}{2}}^{(0)}(t) + 21U_i^{(0)}(t) + 14u_{i+\frac{1}{2}}^{(0)}(t) \right.   \\
		\left.-45U_{i+1}^{(0)}(t) +16 u_{i+\frac{3}{2}}^{(0)}(t) +15U_{i+2}^{(0)}(t) -13u_{i+\frac{5}{2}}^{(0)}(t) + 6U_{i+3}^{(0)}(t) -2 u_{i+\frac{7}{2}}^{(0)}(t)\right) \\
		= \frac{\partial^3}{\partial x^3}u^{(0)}(t,x_{i+\frac{1}{2}}) - \Delta x \frac{\partial^4}{\partial x^4}u^{(0)}(t,x_{i+\frac{1}{2}}) + \mathcal{O}\left(\Delta x^2\right).
\end{multline*}
as well as
\begin{equation}
	\begin{multlined}
		\frac{1}{\Delta x}\left(u_{i-\frac{1}{2}}^{(1)}(t)- 3 U_i^{(1)}(t) +4u_{i+\frac{1}{2}}^{(1)}(t)- 3 U_{i+1}^{(1)}(t) + u_{i+\frac{3}{2}}^{(1)}(t)\right) \\
	= \frac{\Delta x^3}{30}\frac{\partial^4}{\partial x^4}u^{(1)}(t,x_{i+\frac{1}{2}}) + \frac{\Delta x^5}{630}\frac{\partial^6}{\partial x^6}u^{(1)}(t,x_{i+\frac{1}{2}}) +  \mathcal{O}\left(\Delta x^7\right).
	\end{multlined} \label{eq:missing_term}
\end{equation}
Hence, \eqref{eq:3_3_2} is a spatial discretisation of \eqref{eq:th_1} with an error in $\mathcal{O}(\Delta x)$. For the terms appearing in \eqref{eq:3_3_4}, we obtain
\begin{equation*}
	\frac{1}{ \Delta x}\left(\frac{1}{2}\left(u_{i+\frac{1}{2}}^{(0)}(t)\right)^2 -\frac{1}{2}\left(u_{i-\frac{1}{2}}^{(0)}(t)\right)^2  \right)  = \frac{1}{\Delta x}\int_{\mathcal{C}_i} u^{(0)}(t,x) \frac{\partial}{\partial x}u^{(0)}(t,x) \mathrm{d}x
\end{equation*}
and
\begin{equation}
	\label{eq:th_proof_new}
	\begin{multlined}
		\frac{2}{\Delta x^2} \left(-2u_{i-\frac{1}{2}}^{(0)}(t) + 6U_i^{(0)}(t)- 4u_{i+\frac{1}{2}}^{(0)}(t) - 3 U_{i+1}^{(0)}(t) +5 u_{i+\frac{3}{2}}^{(0)}(t) - 3 U_{i+2}^{(0)}(t) + u_{i+\frac{5}{2}}^{(0)(t)}\right) \\
	= \frac{\partial^2}{\partial x^2} u^{(0)}(t,x_{i+\frac{1}{2}}) + \frac{\Delta x}{2}\frac{\partial^3}{\partial x^3}u^{(0)}(t,x_{i+\frac{1}{2}}) + \mathcal{O}\left(\Delta x^2\right),
	\end{multlined}
\end{equation}
which implies
\begin{multline*}
	\frac{2}{\Delta x^3}\left(\bigg(-2u_{i-\frac{1}{2}}^{(0)}(t) + 6U_i^{(0)}(t)- 4u_{i+\frac{1}{2}}^{(0)}(t) - 3 U_{i+1}^{(0)}(t) +5 u_{i+\frac{3}{2}}^{(0)}(t) \right. \\
			\left.- 3 U_{i+2}^{(0)}(t) + u_{i+\frac{5}{2}}^{(0)}(t)\right) - \left(-2u_{i-\frac{3}{2}}^{(0)}(t) + 6U_{i-1}^{(0)}(t)- 4u_{i-\frac{1}{2}}^{(0)}(t) - 3 U_{i}^{(0)} (t) \right.\\
			\left.\left.+5 u_{i+\frac{1}{2}}^{(0)}(t) - 3 U_{i+1}^{(0)}(t) + u_{i+\frac{3}{2}}^{(0)}(t)\right)\right) = \frac{1}{\Delta x}\left(\frac{\partial^2}{\partial x^2} u^{(0)}(t,x_{i+\frac{1}{2}}) - \frac{\partial^2}{\partial x^2} u^{(0)}(t,x_{i-\frac{1}{2}})\right)\\
			 +\frac{1}{2}\left(\frac{\partial^3}{\partial x^3}u^{(0)}(t,x_{i+\frac{1}{2}})-\frac{\partial^3}{\partial x^3}u^{(0)}(t,x_{i-\frac{1}{2}}) \right) + \mathcal{O}(\Delta x).
\end{multline*}
In addition, we have
\begin{equation*}
	\frac{1}{2}\left(\frac{\partial^3}{\partial x^3}u^{(0)}(t,x_{i+\frac{1}{2}})-\frac{\partial^3}{\partial x^3}u^{(0)}(t,x_{i-\frac{1}{2}}) \right) =  \mathcal{O}(\Delta x),
\end{equation*}
where we used the Lipschitz continuity of the derivatives of $u^{(0)}$. Thus, \eqref{eq:3_3_4} is a spatial discretisation of \eqref{eq:th_2} with an error in $\mathcal{O}(\Delta x)$. Equation \eqref{eq:3_3_15} is a discretisation of \eqref{eq:th_3} with an error in $\mathcal{O}(\Delta x^2)$, since 
\begin{multline*}
		\frac{1}{\Delta x}\left(u_{i-\frac{1}{2}}^{(0)}(t) - 3U_i^{(0)}(t) + 3 U_{i+1}^{(0)}(t) - u_{i+\frac{3}{2}}^{(0)}(t)\right) = \frac{\partial}{\partial x}u^{(0)}(t,x_{i+\frac{1}{2}}) - \frac{\Delta x^2}{12} \frac{\partial^3}{\partial x^3}u^{(0)}(t,x_{i+\frac{1}{2}}) + \mathcal{O}(\Delta x^4).
\end{multline*}
Furthermore, 
\begin{align*}
	\frac{1}{\Delta x} \left(u_{i+\frac{1}{2}}^{(0)}(t) - u_{i-\frac{1}{2}}^{(0)}(t) \right) &= \frac{1}{\Delta x} \left(u^{(0)}(t,x_{i+\frac{1}{2}}) - u^{(0)}(t,x_{i-\frac{1}{2}}) \right) = \frac{1}{\Delta x}\int_{\mathcal{C}_i} \frac{\partial}{\partial x}u^{(0)}(t,x) \mathrm{d}x
\end{align*}
holds and it follows that \eqref{eq:3_3_13} is a discretisation of \eqref{eq:th_4} with zero discretisation error. From \eqref{eq:th_proof_new} we immediately obtain that \eqref{eq:3_3_16} is a discretisation of \eqref{eq:th_5} with an error in $\mathcal{O}(\Delta x)$, while 
\begin{multline*}
	\frac{1}{\Delta x^2}\bigg(\left(u_{i-\frac{1}{2}}^{(0)}(t) - 3U_i^{(0)}(t) + 3 U_{i+1}^{(0)}(t) - u_{i+\frac{3}{2}}^{(0)}(t)\right) - \left(u_{i-\frac{3}{2}}^{(0)}(t) - 3U_{i-1}^{(0)}(t) \right.\left.+ 3 U_{i}^{(0)}(t) - u_{i+\frac{1}{2}}^{(0)}(t)\right)\bigg) \\= \frac{1}{\Delta x}\int_{\mathcal{C}_i} \frac{\partial^2}{\partial x^2}u^{(0)}(t,x) \mathrm{d}x + \mathcal{O}(\Delta x^2),
\end{multline*}
gives that \eqref{eq:3_3_14} is a discretisation of \eqref{eq:th_6} with an error in $\mathcal{O}(\Delta x^2)$.
\end{proof}

To sum up the findings: The formal asymptotic analysis indicates that the \texttt{classical} scheme discretised in time with an IMEX Runge-Kutta scheme of type ARS is AP for the component $u$. If, in addition, the IMEX Runge-Kutta scheme is GSA, then it indicates that the scheme is AP. Moreover, the analysis suggests that the truncation error of the limit scheme is in $\mathcal{O}(\Delta x)$. Numerically we observe that the scheme is indeed only first order if the relaxation parameter is  sufficiently small.

\begin{remark}
	The asymptotic analysis yields the additional equation
\begin{align}
	0 &= \frac{1}{2} \left(D^+u_{i+\frac{1}{2}}^{(0)[j]} - D^-u_{i+\frac{1}{2}}^{(0)[j]} \right)_{j=1,\dots,s} \label{eq:unnat}
\end{align}
originating from \eqref{eq:tim_dis_u}.
If the IMEX Runge-Kutta scheme is GSA, this implies
\begin{align*}
	0 &= \frac{1}{2} \left(D^+u_{i+\frac{1}{2}}^{(0)n+1} - D^-u_{i+\frac{1}{2}}^{(0)n+1} \right),
\end{align*}
which is a discretisation of
\begin{align*}
	\frac{\Delta x^3}{30}\frac{\partial^4}{\partial x^4}u^{(0)}(t,x) = 0
\end{align*}
with an error in $\mathcal{O}(\Delta x^5)$. Since this condition does not correspond to any part of the limit equations, we refer to it as an unnatural condition. The formal asymptotic analysis gives no insight into whether and how the unnatural condition affects the AP property of \texttt{classical} scheme.
\end{remark}

\section{Alternative Active Flux schemes}
\label{sec:alternative_scheme}
Numerical experiments conducted for the \texttt{classical} scheme (discretised in time with IMEX Runge-Kutta schemes of type ARS) confirm the AP property, but also indicate that some of the advantageous properties of the scheme vanish as $\varepsilon \rightarrow 0$. In particular, the \texttt{classical} scheme is not AA, and the discrete energy is not conserved well for small $\varepsilon$. This leads us to develop and analyse several alternative Active Flux schemes with different point value update equations. These schemes differ from the \texttt{classical} scheme only in the way the term
	\begin{equation*}
		F_2 \frac{\partial}{\partial x}q(t,x)
	\end{equation*}
  	in \eqref{eq:qua_lin} is discretised, while all other parts of the scheme remain the same. 

\subsection{Optimised version of the classical scheme}
The first approach to develop Active Flux schemes with improved properties is to optimise the \texttt{classical} scheme by modifying the matrices $F_2^\pm$. Observe that the terms of factor $\frac{1}{\varepsilon}$ in \eqref{eq:jp_and_jm_2} give rise to expressions of the form \eqref{eq:dp_m_dm} which is a discretisation of 
\begin{align*}
	\frac{\Delta x^3}{30}\frac{\partial^4}{\partial x^4}v(t,x_{i+\frac{1}{2}}) + \mathcal{O}(\Delta x^5).
\end{align*}
(see \eqref{eq:missing_term}).These appear in \eqref{eq:tim_dis_u}--\eqref{eq:tim_dis_u2} and \eqref{eq:tim_dis_pv_psi}--\eqref{eq:tim_dis_pv_psi2}. If we replace the two $\frac{1}{\varepsilon}$ elements on the diagonal of $F_2^\pm$ with zeros, \eqref{eq:pv_cll} is still a discretisation of \eqref{eq:qua_lin}. Numerical experiments show that the resulting scheme is at most second order accurate and performs worse than the \texttt{classical} scheme for $\varepsilon$ close to $1$. However, for sufficiently small $\varepsilon$, this scheme is slightly more accurate than the \texttt{classical} scheme. To refine this approach, we only replace the first element in the first row (the ``1-1-element'') of the matrices $F_2^\pm$ by zero, since it is the one responsible for the unnatural condition \eqref{eq:unnat}. The 3-3-element of $F_2^\pm$ scales with $\frac{1}{\varepsilon}$, while the 3-1-element scales with $\frac{1}{\varepsilon^2}$. As a result, if $\varepsilon \rightarrow 0$, the relation between the terms in update equation of the point values changes. By changing the scaling of the 3-3-element to $\frac{1}{\varepsilon^2}$, the relation of the terms in the update equation becomes independent of $\varepsilon$. These arguments will become clearer below when the asymptotic analysis is performed. Following this new approach, the splitting of the Jacobian now reads
\begin{align}
	F_{2,\mathrm{optimised}}^+ = \begin{pmatrix}
		0 & 0 & \frac{1}{2} \\
		0 & 0 & 0 \\
		\frac{1}{2\varepsilon^2} & 0 & \frac{1}{2\varepsilon^2}
	\end{pmatrix} && \text{and} && F_{2,\mathrm{optimised}}^- = \begin{pmatrix}
		0 & 0 & \frac{1}{2} \\
		0 & -\frac{1}{\varepsilon^2} & 0 \\
		\frac{1}{2\varepsilon^2} & 0 & -\frac{1}{2\varepsilon^2}
	\end{pmatrix}.\label{eq:f2optimised}
\end{align}
Observe that the eigenvalues of $F_{2,\mathrm{optimised}}^\pm$ no longer have a definite sign. However, ad hoc modifications of the upwinding matrices are very successful in restoring the correct asymptotics while retaining stability in other situations already, e.g.\ as ``low Mach fixes'' in the context of the incompressible limit of the Euler equations, see e.g.\ \cite{li_gu_2013}. We call the resulting scheme the Active Flux scheme with optimised individual upwinding or short the \texttt{optimised} scheme. 

Performing a formal asymptotic analysis for the \texttt{optimised} scheme, discretised in time with an IMEX Runge-Kutta scheme of type ARS, allows us to obtain \eqref{eq:3_3_3}, \eqref{eq:3_3_5}, and \eqref{eq:3_3_6} for the terms of order $0$, since the update equation of the cell averages remains unchanged. In addition, we get \eqref{eq:3_3_7} (since \eqref{eq:tim_dis_pv_p}--\eqref{eq:tim_dis_pv_p2} and thus \eqref{eq:av_faa_3_1}--\eqref{eq:av_faa_3_2} remain valid) as well as
\begin{subequations}
    \label{eq:4_3_1}
	\begin{align}
		u_{i+\frac{1}{2}}^{[k](0)} &= 
	\begin{aligned}[t]
		&u^n_{i+\frac{1}{2}} + \Delta t \sum_{j=1}^{k}\tilde a_{k+1,j} \left(-\max\left((u_{i+\frac{1}{2}}^{(0)[j-1]},0\right)D^+u_{i+\frac{1}{2}}^{(0)[j-1]}\right.\\
        &\left.-\min\left(u_{i+\frac{1}{2}}^{(0)[j-1]},0\right)D^-u_{i+\frac{1}{2}}^{(0)[j-1]} \right) +\Delta t \sum_{j=1}^k a_{k,j}\left( -  \frac{1}{2}\left(D^+\psi_{i+\frac{1}{2}}^{(0)[j]} + D^-\psi_{i+\frac{1}{2}}^{(0)[j]} \right) \right) ,
	\end{aligned}\\
	u_{i+\frac{1}{2}}^{n+1(0)} &= \begin{aligned}[t]
		&u^n_{i+\frac{1}{2}} + \Delta t \sum_{j=1}^{\sigma}\tilde b_j \left(-\max\left(u_{i+\frac{1}{2}}^{(0)[j-1]},0\right)D^+u_{i+\frac{1}{2}}^{(0)[j-1]}\right.\\
        &\left.-\min\left(u_{i+\frac{1}{2}}^{(0)[j-1]},0\right)D^-u_{i+\frac{1}{2}}^{(0)[j-1]} \right) +\Delta t \sum_{j=1}^s b_{j}\left( -  \frac{1}{2}\left(D^+\psi_{i+\frac{1}{2}}^{(0)[j]} + D^-\psi_{i+\frac{1}{2}}^{(0)[j]} \right) \right) 
	\end{aligned}
	\end{align}
\end{subequations}
and 
\begin{align}
	\rho_{i+\frac{1}{2}}^{(0)[j]} =\frac{1}{2} \left(D^+u_{i+\frac{1}{2}}^{(0)[j]} + D^-u_{i+\frac{1}{2}}^{(0)[j]} \right) + \frac{1}{2} \left(D^+\psi_{i+\frac{1}{2}}^{(0)[j]} - D^-\psi_{i+\frac{1}{2}}^{(0)[j]} \right) 
	\label{eq:4_3_2}
\end{align}
instead of \eqref{eq:3_3_1} and \eqref{eq:3_3_8}, respectively. The absence of the $u^{(1)}$-terms in \eqref{eq:4_3_1} and the presence of the $\psi^{(0)}$-terms in \eqref{eq:4_3_2} are due to our modification. The unnatural condition \eqref{eq:unnat} is no longer present. It follows in a similar way as before that \eqref{eq:4_3_1} is a temporal discretisation (using an type ARS IMEX Runge-Kutta scheme) of 
\begin{equation}
	\begin{aligned}
		\frac{\mathrm{d}}{\mathrm{d}t} u_{i+\frac{1}{2}}^{(0)}(t) = \begin{aligned}[t]
			& -\bigg(\max\left(u_{i+\frac{1}{2}}^{(0)}(t) ,0\right) \frac{1}{\Delta x}\left(2u_{i-\frac{1}{2}}^{(0)}(t) - 6U_i^{(0)}(t) + 4 u_{i+\frac{1}{2}}^{(0)}(t)\right) \\
			& + \min\left(u_{i+\frac{1}{2}}^{(0)}(t) ,0\right) \frac{1}{\Delta x}\left(-4u_{i+\frac{1}{2}}^{(0)}(t) + 6 U_{i+1}^{(0)}(t) - 2 u_{i+\frac{3}{2}}^{(0)}(t)  \right) \bigg)  \\
			& -\frac{1}{\Delta x^3}\bigg(-u_{i-\frac{3}{2}}^{(0)}(t) + 3U_{i-1}^{(0)}(t) -14u_{i-\frac{1}{2}}^{(0)}(t) + 21U_i^{(0)}(t) + 14u_{i+\frac{1}{2}}^{(0)}(t) \\
			& -45U_{i+1}^{(0)}(t)  +16 u_{i+\frac{3}{2}}^{(0)}(t) +15U_{i+2}^{(0)}(t) -13u_{i+\frac{5}{2}}^{(0)}(t) + 6U_{i+3}^{(0)}(t) -2 u_{i+\frac{7}{2}}^{(0)}(t)\bigg) \\
			& - \frac{1}{\Delta x^2} \left(-\Lambda_{i-\frac{1}{2}}^{(0)}(t) -8 \Lambda_{i+\frac{1}{2}}^{(0)}(t) +7 \Lambda_{i+\frac{3}{2}}^{(0)}(t) + 2\Lambda_{i+\frac{5}{2}}^{(0)}(t)  \right)
		\end{aligned}
	\end{aligned}
	\label{eq:4_3_3}
\end{equation}
and \eqref{eq:3_3_3} is a discretisation of
\begin{equation}
	\begin{aligned}
	\frac{\mathrm{d}}{\mathrm{d}t}U_i^{(0)}(t) = \begin{aligned}[t]
			& -\frac{1}{\Delta x}\left(\frac{1}{2}\left(u_{i+\frac{1}{2}}^{(0)}(t)\right)^2 -\frac{1}{2}\left(u_{i-\frac{1}{2}}^{(0)}(t)\right)^2  \right)  \\
			&- \frac{2}{\Delta x^3}\left(\bigg(-2u_{i-\frac{1}{2}}^{(0)}(t) + 6U_i^{(0)}(t)- 4u_{i+\frac{1}{2}}^{(0)}(t) - 3 U_{i+1}^{(0)}(t) +5 u_{i+\frac{3}{2}}^{(0)}(t) \right. \\
			&\left.- 3 U_{i+2}^{(0)}(t) + u_{i+\frac{5}{2}}^{(0)}(t)\right) - \left(-2u_{i-\frac{3}{2}}^{(0)}(t) + 6U_{i-1}^{(0)}(t)- 4u_{i-\frac{1}{2}}^{(0)}(t) - 3 U_{i}^{(0)} (t) \right.\\
			&\left.\left.+5 u_{i+\frac{1}{2}}^{(0)}(t) - 3 U_{i+1}^{(0)}(t) + u_{i+\frac{3}{2}}^{(0)}(t)\right)\right) \\
			&- \frac{1}{\Delta x^2}\left(\left(-4 \Lambda_{i+\frac{1}{2}}^{(0)}(t) - 2 \Lambda_{i+\frac{3}{2}}^{(0)}(t)\right)  - \left(-4 \Lambda_{i-\frac{1}{2}}^{(0)}(t) - 2 \Lambda_{i+\frac{1}{2}}^{(0)}(t)\right)\right),
		\end{aligned}
	\end{aligned}
	\label{eq:4_3_4}
\end{equation}
where 
\begin{align*}
	\Lambda_{i+\frac{1}{2}}^{(0)}(t) = \frac{1}{2}\left(D^+\psi_{i+\frac{1}{2}}^{(0)}(t)  - D^-\psi_{i+\frac{1}{2}}^{(0)}(t)\right).
\end{align*}
If the time discretisation scheme is GSA, then we additionally obtain the semi-discretisations \eqref{eq:3_3_13} and, instead of \eqref{eq:3_3_14}
\begin{align}
	\begin{aligned}
		\Psi_i^{(0)}(t) =\begin{aligned}[t]
			&\frac{1}{\Delta x^2} \left(\left(u_{i-\frac{1}{2}}^{(0)}(t) - 3U_i^{(0)}(t) + 3U_{i+1}^{(0)}(t) - u_{i+\frac{3}{2}}^{(0)}(t)  \right) \right. \\
		 &\left.- \left(u_{i-\frac{3}{2}}^{(0)}(t) - 3U_{i-1}^{(0)}(t) + 3U_{i}^{(0)}(t) - u_{i+\frac{1}{2}}^{(0)}(t)  \right)  \right) \\
		 &+ \frac{1}{\Delta x^2} \left(\left( \psi_{i-\frac{1}{2}}^{(0)}(t) - 3 \Psi_i^{(0)}(t) + 4 \psi_{i+\frac{1}{2}}^{(0)}(t) - 3 \Psi_{i+1}^{(0)}(t) + \psi_{i+\frac{3}{2}}^{(0)}(t) \right)\right.\\
		 &-\left. \left( \psi_{i-\frac{3}{2}}^{(0)}(t) - 3 \Psi_{i-1}^{(0)}(t) + 4 \psi_{i-\frac{1}{2}}^{(0)}(t) - 3 \Psi_{i}^{(0)}(t) + \psi_{i+\frac{1}{2}}^{(0)}(t) \right)\right)
		\end{aligned}
	\end{aligned}
	\label{eq:4_3_8}
\end{align}
as well as 
\begin{equation}
	\begin{aligned}
		\rho_{i+\frac{1}{2}}^{(0)}(t) =\begin{aligned}[t]
		&\frac{1}{\Delta x}\left(u_{i-\frac{1}{2}}^{(0)}(t) - 3U_i^{(0)}(t) + 3 U_{i+1}^{(0)}(t) - u_{i+\frac{3}{2}}^{(0)}(t)\right) \\
		&+\frac{1}{\Delta x}\left( \psi_{i-\frac{1}{2}}^{(0)}(t) - 3 \Psi_i^{(0)}(t) + 4 \psi_{i+\frac{1}{2}}^{(0)}(t) - 3 \Psi_{i+1}^{(0)}(t) + \psi_{i+\frac{3}{2}}^{(0)}(t) \right)
	\end{aligned} 
	\end{aligned}
	\label{eq:4_3_9}
\end{equation}
and
\begin{equation}
	\begin{aligned}
		\psi_{i+\frac{1}{2}}^{(0)}(t) = \begin{aligned}[t]
				&\frac{2}{\Delta x^2} \left(-2u_{i-\frac{1}{2}}^{(0)}(t) + 6U_i^{(0)}(t)- 4u_{i+\frac{1}{2}}^{(0)}(t) - 3 U_{i+1}^{(0)}(t) +5 u_{i+\frac{3}{2}}^{(0)}(t) \right.  \\
				& \left. - 3 U_{i+2}^{(0)}(t) + u_{i+\frac{5}{2}}^{(0)}(t)\right) + \frac{1}{\Delta x} \left(-4\Lambda_{i+\frac{1}{2}}^{(0)}(t) -2 \Lambda_{i+\frac{3}{2}}^{(0)}(t) \right)
			\end{aligned}
	\end{aligned}
	\label{eq:4_3_10}
\end{equation}
instead of \eqref{eq:3_3_15} and \eqref{eq:3_3_16}, respectively.
Assuming $u^{(0)},\psi^{(0)} :\R^+_0 \times \R \rightarrow \R$  are regular enough and interpolate the point values and cell averages, the additional $\Lambda$-terms are in $\mathcal{O}(\Delta x^3)$. The difference between \eqref{eq:4_3_1} and \eqref{eq:3_3_1} is the term \eqref{eq:missing_term}. This term is in $\mathcal{O}(\Delta x^3)$. Therefore, a similar result as Theorem \ref{th:3_3_1} holds and the formal asymptotic analysis indicates that the \texttt{optimised} scheme is AP for the $u$ component. If the IMEX Runge-Kutta scheme is GSA, then the scheme is AP. 

For $\varepsilon$ close to $1$ the \texttt{optimised} scheme performs similarly as the \texttt{classical} scheme, as will be seen in the experiments of Section \ref{sec:num_exp}. However, if $\varepsilon \rightarrow 0$, the properties of the \texttt{optimised} scheme do not deteriorate. In particular, the scheme achieves third-order accuracy independently of $\varepsilon$, i.e.\ it is AA, although the asymptotic analysis continues to show presence of terms that are only $\mathcal O(\Delta x)$ approximations (e.g. in \eqref{eq:4_3_3}).

\subsection{Central Active Flux schemes }

For time-explicit schemes, upwinding is typically used to ensure the stability of the scheme. It is common to use a central discretisation for the term discretised in time implicitly, i.e.\ in our case we would be thus led to choosing $F_2^\pm = \frac{1}{2}F_2$. The resulting matrices read
\begin{align}
	\label{eq:jp_and_jm_2_central}
	F_{2,\mathrm{central}}^+ = \begin{pmatrix}
		0 & 0 & \frac{1}{2} \\
		0 & -\frac{1}{2\varepsilon^2} & 0 \\
		\frac{1}{2 \varepsilon^2} & 0 & 0
	\end{pmatrix} &&\text{and} && F_{2,\mathrm{central}}^- = \begin{pmatrix}
		0 & 0 & \frac{1}{2} \\
		0 & -\frac{1}{2\varepsilon^2} & 0 \\
		\frac{1}{2 \varepsilon^2} & 0 & 0
	\end{pmatrix}.
\end{align}
We refer to an Active Flux scheme with \eqref{eq:jp_and_jm_2_central} instead of \eqref{eq:jp_and_jm_2} in \eqref{eq:pv_cll} as the central Active Flux scheme or short the \texttt{central} scheme. A formal asymptotic analysis for the time-discrete \texttt{central} scheme yields, besides \eqref{eq:3_3_3}, \eqref{eq:3_3_5}, \eqref{eq:3_3_6}, and \eqref{eq:3_3_8}, also 
\begin{subequations}
     \label{eq:4_2_1}
	\begin{align}
		u_{i+\frac{1}{2}}^{[k](0)} &= 
	\begin{aligned}[t]
		&u^n_{i+\frac{1}{2}} + \Delta t \sum_{j=1}^{k}\tilde a_{k+1,j} \left(-\max\left((u_{i+\frac{1}{2}}^{(0)[j-1]},0\right)D^+u_{i+\frac{1}{2}}^{(0)[j-1]}\right.\\
        &\left.-\min\left(u_{i+\frac{1}{2}}^{(0)[j-1]},0\right)D^-u_{i+\frac{1}{2}}^{(0)[j-1]} \right) \\
		&+\Delta t \sum_{j=1}^k a_{k,j}\left( -  \frac{1}{2}\left(D^+\psi_{i+\frac{1}{2}}^{(0)[j]} + D^-\psi_{i+\frac{1}{2}}^{(0)[j]} \right) \right), 
	\end{aligned}\\
	u_{i+\frac{1}{2}}^{n+1(0)} &= \begin{aligned}[t]
		&u^n_{i+\frac{1}{2}} + \Delta t \sum_{j=1}^{\sigma}\tilde b_j \left(-\max\left(u_{i+\frac{1}{2}}^{(0)[j-1]},0\right)D^+u_{i+\frac{1}{2}}^{(0)[j-1]}\right.\\
        &\left.-\min\left(u_{i+\frac{1}{2}}^{(0)[j-1]},0\right)D^-u_{i+\frac{1}{2}}^{(0)[j-1]} \right) \\
		&+\Delta t \sum_{j=1}^s b_{j}\left( -  \frac{1}{2}\left(D^+\psi_{i+\frac{1}{2}}^{(0)[j]} + D^-\psi_{i+\frac{1}{2}}^{(0)[j]} \right) \right) 
	\end{aligned}
	\end{align}
\end{subequations}
instead of \eqref{eq:3_3_1} and
\begin{align}
	\psi_{i+\frac{1}{2}}^{(0)[j]} = \frac{1}{2} \left(D^+\rho_{i+\frac{1}{2}}^{(0)[j]} + D^-\rho_{i+\frac{1}{2}}^{(0)[j]} \right). \label{eq:4_2_2}
\end{align}
instead of \eqref{eq:3_3_7}.
In particular, the scheme does not give rise to any unnatural conditions (like \eqref{eq:unnat}) for the terms of order $0$. We can recover the following time-continuous discretisations:
\begin{equation}
	\begin{aligned}
		\frac{\mathrm{d}}{\mathrm{d}t} u_{i+\frac{1}{2}}^{(0)}(t) = \begin{aligned}[t]
			& -\left(\max\left(u_{i+\frac{1}{2}}^{(0)}(t) ,0\right) \frac{1}{\Delta x}\left(2u_{i-\frac{1}{2}}^{(0)}(t) - 6U_i^{(0)}(t) + 4 u_{i+\frac{1}{2}}^{(0)}(t)\right)\right. \\
			& + \left.\min\left(u_{i+\frac{1}{2}}^{(0)}(t) ,0\right) \frac{1}{\Delta x}\left(-4u_{i+\frac{1}{2}}^{(0)}(t) + 6 U_{i+1}^{(0)}(t) - 2 u_{i+\frac{3}{2}}^{(0)}(t)  \right) \right)  \\
			&-\frac{1}{\Delta x^3}\left(u_{i-\frac{5}{2}}^{(0)}(t) -3U_{i-2}^{(0)}(t) + 6 u_{i-\frac{3}{2}}^{(0)}(t) -6U_{i-1}^{(0)}(t)\right. \\
			&- 15u_{i-\frac{1}{2}}^{(0)}(t) + 33 U_i^{(0)}(t) - 33U_{i+1}^{(0)}(t) + 15u_{i+\frac{3}{2}}^{(0)}(t)\\
			& \left.+ 6U_{i-1}^{(0)}(t) - 6u_{i+\frac{5}{2}}^{(0)}(t) + 3U_{i+3}^{(0)}(t) - u_{i-\frac{7}{2}}^{(0)}(t)  \right)	
	\end{aligned}
	\end{aligned}
	\label{eq:4_2_3}
\end{equation} 
instead of \eqref{eq:3_3_2},
\begin{equation}
	\begin{aligned}
	\frac{\mathrm{d}}{\mathrm{d}t}U_i^{(0)}(t) = \begin{aligned}[t]
			& -\frac{1}{2 \Delta x}\left(\left(u_{i+\frac{1}{2}}^{(0)}(t)\right)^2 -\left(u_{i-\frac{1}{2}}^{(0)}(t)\right)^2  \right) - \frac{1}{\Delta x^3}\left(\left(u_{i-\frac{3}{2}}^{(0)}(t)  - 3U_{i-1}^{(0)}(t) \right. \right. \\
			& +3 u_{i-\frac{1}{2}}^{(0)}(t) +3U_i^{(0)}(t) -8 u_{i+\frac{1}{2}}^{(0)}(t) +3U_{i+1}^{(0)}(t) +3u_{i+\frac{3}{2}}^{(0)}(t)\\
			&\left.  - 3U_{i+2}^{(0)}(t)  + u_{i+\frac{5}{2}}^{(0)}(t)  \right) - \left(u_{i-\frac{5}{2}}^{(0)}(t)  - 3U_{i-2}^{(0)}(t) +3 u_{i-\frac{3}{2}}^{(0)}(t) \right. \\
			&\left. \left. +3U_{i-1}^{(0)}(t) -8 u_{i-\frac{1}{2}}^{(0)}(t) +3U_{i}^{(0)}(t) +3u_{i+\frac{1}{2}}^{(0)}(t)  - 3U_{i+1}^{(0)}(t)  + u_{i+\frac{3}{2}}^{(0)}(t)  \right)\right).
		\end{aligned}
	\end{aligned}
	\label{eq:4_2_4}
\end{equation}
instead of \eqref{eq:3_3_4}, and
\begin{equation}
	\begin{aligned}
		\psi_{i+\frac{1}{2}}^{(0)}(t) = \begin{aligned}[t]
 			&\frac{1}{\Delta x^2}\left(u_{i-\frac{3}{2}}^{(0)}(t)  - 3U_{i-1}^{(0)}(t) +3 u_{i-\frac{1}{2}}^{(0)}(t) +3U_i^{(0)}(t) -8 u_{i+\frac{1}{2}}^{(0)}(t)\right. \\
 			&\left.+3U_{i+1}^{(0)}(t) +3u_{i+\frac{3}{2}}^{(0)[j]}  - 3U_{i+2}^{(0)}(t)  + u_{i+\frac{5}{2}}^{(0)}(t)  \right).
 			\end{aligned}
	\end{aligned}
	\label{eq:4_2_7}
\end{equation}
instead of \eqref{eq:3_3_16} as well as \eqref{eq:3_3_15}, \eqref{eq:3_3_13}, and \eqref{eq:3_3_14}. As before, recovering the latter four requires the assumption that the IMEX Runge-Kutta scheme is GSA.
Observe the new, more symmetric discretisations of $\frac{\partial^3}{\partial x^3}u$ and $\frac{\partial^2}{\partial x^2}u$ in \eqref{eq:4_2_4} and \eqref{eq:4_2_7}, respectively.
 Similarly to the proof of Theorem \ref{th:3_3_1}, one can show that \eqref{eq:4_2_3}, \eqref{eq:3_3_15}, and \eqref{eq:4_2_7} are a discretisation of \eqref{eq:th_1},  \eqref{eq:th_3}, and \eqref{eq:th_5} with an error in $\mathcal{O}(\Delta x^2)$ and \eqref{eq:4_2_4}, \eqref{eq:3_3_13}, and \eqref{eq:3_3_14} are a discretisation of \eqref{eq:th_2},  \eqref{eq:th_4}, and \eqref{eq:th_6} with the same discretisation error. Hence, the \texttt{central} scheme is AP for the $u$ component and AP for all components if the IMEX Runge-Kutta scheme is GSA. Numerically this property can be verified. We find that the \texttt{central} scheme is second order accurate for both large and small $\varepsilon$, i.e.\ the scheme is AA since it retains its order of accuracy, but it does not reach the usual third order of Active Flux methods. The present analysis does not offer a conclusive theoretical explanation for this loss of accuracy.

It is possible to modify the central scheme to obtain a third-order AA scheme by considering the following alternatives to $D^\pm$ inspired by \cite{duan_2026}:
\begin{subequations}
	\label{eq:mod_dp_dm}
	\begin{align}
		D^+_\mathrm{mod}v_{i+\frac{1}{2}}(t) &= \frac{1}{6 \Delta x}\left(5v_{i-\frac{1}{2}}(t) - 18 V_i(t) + 8 v_{i+\frac{1}{2}}(t) + 6V_{i+1}(t) - v_{i+\frac{3}{2}}(t) \right) \\
		D^-_\mathrm{mod}v_{i+\frac{1}{2}}(t) &= \frac{1}{6 \Delta x}\left(v_{i-\frac{1}{2}}(t) - 6 V_i(t) - 8 v_{i+\frac{1}{2}}(t) + 18V_{i+1}(t) -5 v_{i+\frac{3}{2}}(t) \right).
	\end{align}
\end{subequations}
We replace $D^\pm $ by $D^\pm_\mathrm{mod}$ only in the part that is discretised in time implicitly.
We name the scheme in which we replace all $D^\pm$ in the implicit part with the new formulas \eqref{eq:mod_dp_dm} the first modified central Active Flux scheme or \texttt{mod1-central} scheme. In fact, we find it sufficient to only make the replacement in the update equations of the point values of $\rho$ and $\psi$. This scheme is referred to as the second modified central Active Flux scheme or \texttt{mod2-central} scheme.

Since the two modified central schemes are closely related, the formal asymptotic analysis of the schemes is almost identical. Following the usual approach, we obtain \eqref{eq:3_3_3}, \eqref{eq:3_3_5}, \eqref{eq:3_3_6}, and 
\begin{subequations}
\label{eq:4_4_1}
\begin{align}
		u_{i+\frac{1}{2}}^{[k](0)} &= 
	\begin{aligned}[t]
		&u^n_{i+\frac{1}{2}} + \Delta t \sum_{j=1}^{k}\tilde a_{k+1,j} \left(-\max\left((u_{i+\frac{1}{2}}^{(0)[j-1]},0\right)D^+u_{i+\frac{1}{2}}^{(0)[j-1]}\right.\\
        &\left.
		-\min\left(u_{i+\frac{1}{2}}^{(0)[j-1]},0\right)D^-u_{i+\frac{1}{2}}^{(0)[j-1]} \right) \\
		&+\Delta t \sum_{j=1}^k a_{k,j}\left( -  \frac{1}{2}\left(D^+_\mathrm{mod}\psi_{i+\frac{1}{2}}^{[j](0)} + D^-_\mathrm{mod}\psi_{i+\frac{1}{2}}^{[j](0)} \right) \right) ,
	\end{aligned}\\
	u_{i+\frac{1}{2}}^{n+1(0)} &= \begin{aligned}[t]
		&u^n_{i+\frac{1}{2}} + \Delta t \sum_{j=1}^{\sigma}\tilde b_j \left(-\max\left(u_{i+\frac{1}{2}}^{[j-1](0)},0\right)D^+u_{i+\frac{1}{2}}^{[j-1](0)}\right.\\
        &\left.
		-\min\left(u_{i+\frac{1}{2}}^{[j-1](0)},0\right)D^-u_{i+\frac{1}{2}}^{[j-1](0)} \right) \\
		&+\Delta t \sum_{j=1}^s b_{j}\left( -  \frac{1}{2}\left(D^+_\mathrm{mod}\psi_{i+\frac{1}{2}}^{[j](0)} + D^-_\mathrm{mod}\psi_{i+\frac{1}{2}}^{[j](0)} \right) \right).
	\end{aligned}
	\end{align} 
\end{subequations}
as well as 
\begin{align}
	\rho_{i+\frac{1}{2}}^{(0)[j]} = \frac{1}{2}\left( D^+_\mathrm{mod}u^{(0)[j]}_{i+\frac{1}{2}} + D^-_\mathrm{mod}u^{(0)[j]}_{i+\frac{1}{2}}\right)
	\label{eq:4_4_2}
\end{align}
and
\begin{align}
	\psi_{i+\frac{1}{2}}^{(0)[j]} = \frac{1}{2} \left(D^+_\mathrm{mod}\rho_{i+\frac{1}{2}}^{(0)[j]} + D^-_\mathrm{mod}\rho_{i+\frac{1}{2}}^{(0)[j]} \right).
	\label{eq:4_4_3}
\end{align}
for the terms of order $0$ of the \texttt{mod1-central} scheme. For the \texttt{mod2-central} scheme, the analysis yields the same equations with the exception that we get \eqref{eq:4_2_1} instead of \eqref{eq:4_4_1}. For both schemes unnatural conditions of the kind of \eqref{eq:unnat} do not appear.

The equations for the terms of order $0$ of the \texttt{mod1-central} scheme allow us to obtain the following semi-discretisations: 
\begin{equation}
	\begin{aligned}
		\frac{\mathrm{d}}{\mathrm{d}t} u_{i+\frac{1}{2}}^{(0)}(t) = \begin{aligned}[t]
			& -\bigg(\max\left(u_{i+\frac{1}{2}}^{(0)}(t) ,0\right) \frac{1}{\Delta x}\left(2u_{i-\frac{1}{2}}^{(0)}(t) - 6U_i^{(0)}(t) + 4 u_{i+\frac{1}{2}}^{(0)}(t)\right) \\
			& + \min\left(u_{i+\frac{1}{2}}^{(0)}(t) ,0\right) \frac{1}{\Delta x}\left(-4u_{i+\frac{1}{2}}^{(0)}(t) + 6 U_{i+1}^{(0)}(t) - 2 u_{i+\frac{3}{2}}^{(0)}(t)  \right) \bigg)  \\
			&-\frac{1}{8\Delta x^3}\bigg(u_{i-\frac{5}{2}}^{(0)}(t) -4U_{i-2}^{(0)}(t) + 16 u_{i-\frac{3}{2}}^{(0)}(t) -28U_{i-1}^{(0)}(t) \\
			&- 35u_{i-\frac{1}{2}}^{(0)}(t) + 104 U_i^{(0)}(t)  - 104U_{i+1}^{(0)}(t) + 35u_{i+\frac{3}{2}}^{(0)}(t)  \\
			&+ 28U_{i+2}^{(0)}(t) - 16u_{i+\frac{5}{2}}^{(0)}(t) + 4U_{i+3}^{(0)}(t) - u_{i-\frac{7}{2}}^{(0)}(t)  \bigg).		
	\end{aligned}
	\end{aligned}
	\label{eq:4_4_4}
\end{equation}
and
\begin{equation}
	\begin{aligned}
	\frac{\mathrm{d}}{\mathrm{d}t}U_i^{(0)}(t) = \begin{aligned}[t]
			& -\frac{1}{\Delta x}\left(\frac{1}{2}\left(u_{i+\frac{1}{2}}^{(0)}(t)\right)^2 -\frac{1}{2}\left(u_{i-\frac{1}{2}}^{(0)}(t)\right)^2  \right)  \\
			&- \frac{1}{4\Delta x^3}\left(\left(u_{i-\frac{3}{2}}^{(0)}(t)  - 4U_{i-1}^{(0)}(t) + 8 u_{i-\frac{1}{2}}^{(0)} +4U_i^{(0)}(t) -18 u_{i+\frac{1}{2}}^{(0)}(t) +4U_{i+1}^{(0)}(t)  \right. \right. \\
			&\left.+ 8u_{i+\frac{3}{2}}^{(0)}(t) - 4U_{i+2}^{(0)}(t)  + u_{i+\frac{5}{2}}^{(0)}(t)   \right) - \left(u_{i-\frac{5}{2}}^{(0)}(t)  - 4U_{i-2}^{(0)}(t) + 8 u_{i-\frac{3}{2}}^{(0)}\right. \\
			&\left. \left. +4U_{-1}^{(0)}(t) -18 u_{i-\frac{1}{2}}^{(0)}(t) +4U_{i}^{(0)}(t) + 8u_{i+\frac{1}{2}}^{(0)}(t) - 4U_{i+1}^{(0)}(t)  + u_{i+\frac{3}{2}}^{(0)}(t)  \right)\right).
		\end{aligned}
	\end{aligned}
	\label{eq:4_4_5}
\end{equation}
and under the assumption that the time discretisation is GSA \eqref{eq:3_3_13} and
\begin{equation}
	\rho_{i+\frac{1}{2}}^{(0)}(t) = \frac{1}{2\Delta x}\left(u_{i-\frac{1}{2}}^{(0)}(t) - 4U_i^{(0)}(t)) + 4 U_{i+1}^{(0)}(t) - u_{i+\frac{3}{2}}^{(0)}(t) \right),
	\label{eq:4_4_10}
\end{equation}
\begin{equation}
	\begin{aligned}
		\psi_{i+\frac{1}{2}}^{(0)}(t) = \begin{aligned}[t]
 			&\frac{1}{4 \Delta x^2}\left(u_{i-\frac{3}{2}}^{(0)}(t)  - 4U_{i-1}^{(0)}(t) +8 u_{i-\frac{1}{2}}^{(0)}(t) +4U_i^{(0)}(t) -18 u_{i+\frac{1}{2}}^{(0)}(t)\right. \\
 			&\left.+4U_{i+1}^{(0)}(t) +8u_{i+\frac{3}{2}}^{(0)}(t)  - 4U_{i+2}^{(0)}(t)  + u_{i+\frac{5}{2}}^{(0)}(t)  \right)
 			\end{aligned}
	\end{aligned}
	\label{eq:4_4_11}
\end{equation}
and 
\begin{equation}
	\begin{aligned}
		\Psi^{(0)}_i(t) = \begin{aligned}[t]
		&\frac{1}{\Delta 2x^2} \left(\left(u_{i-\frac{1}{2}}^{(0)}(t) - 4U_i^{(0)}(t) + 4U_{i+1}^{(0)}(t) - u_{i+\frac{3}{2}}^{(0)}(t)  \right) \right. \\ 
		 &\left.- \left(u_{i-\frac{3}{2}}^{(0)}(t) - 4U_{i-1}^{(0)}(t) + 4U_{i}^{(0)}(t) - u_{i+\frac{1}{2}}^{(0)}(t)  \right)  \right).
	\end{aligned} 
	\end{aligned}
	\label{eq:4_4_12}
\end{equation}
In case of the \texttt{mod2-central} scheme, we get
\begin{equation}
	\begin{aligned}
		\frac{\mathrm{d}}{\mathrm{d}t} u_{i+\frac{1}{2}}^{(0)}(t) = \begin{aligned}[t]
			& -\bigg(\max\left(u_{i+\frac{1}{2}}^{(0)}(t) ,0\right) \frac{1}{\Delta x}\left(2u_{i-\frac{1}{2}}^{(0)}(t) - 6U_i^{(0)}(t) + 4 u_{i+\frac{1}{2}}^{(0)}(t)\right) \\
			& + \min\left(u_{i+\frac{1}{2}}^{(0)}(t) ,0\right) \frac{1}{\Delta x}\left(-4u_{i+\frac{1}{2}}^{(0)}(t) + 6 U_{i+1}^{(0)}(t) - 2 u_{i+\frac{3}{2}}^{(0)}(t)  \right) \bigg)  \\
			&-\frac{1}{4\Delta x^3}\bigg(u_{i-\frac{5}{2}}^{(0)}(t) -4U_{i-2}^{(0)}(t) + 14 u_{i-\frac{3}{2}}^{(0)}(t) -20 U_{i-1}^{(0)}(t) \\
			&- 31u_{i-\frac{1}{2}}^{(0)}(t) + 80 U_i^{(0)}(t)  - 80U_{i+1}^{(0)}(t) + 31u_{i+\frac{3}{2}}^{(0)}(t)  \\
			&+ 20U_{i+2}^{(0)}(t) - 14u_{i+\frac{5}{2}}^{(0)}(t) + 4U_{i+3}^{(0)}(t) - u_{i-\frac{7}{2}}^{(0)}(t)  \bigg),	
	\end{aligned}
	\end{aligned}
	\label{eq:mod2_4_4_4}
\end{equation}
as well as \eqref{eq:4_4_5}, \eqref{eq:3_3_13}, \eqref{eq:4_4_10}, \eqref{eq:4_4_11}, and \eqref{eq:4_4_12}. Equation \eqref{eq:4_4_4} is a discretisation of \eqref{eq:th_1} with an error in $\mathcal{O}(\Delta x^2)$. The error is due to the discretisation of the nonlinear part, the discretisation error of the third derivative is in $\mathcal{O}(\Delta x^4)$. In addition, \eqref{eq:4_4_10}, \eqref{eq:4_4_11}, \eqref{eq:4_4_5}, and \eqref{eq:4_4_12} are discretisation of \eqref{eq:th_3}, \eqref{eq:th_5}, \eqref{eq:th_2}, and \eqref{eq:th_6} respectively with an error in $\mathcal{O}(\Delta x^4)$. As \eqref{eq:4_4_4}, \eqref{eq:mod2_4_4_4} is also a discretisation of \eqref{eq:th_1} with an error in $\mathcal{O}(\Delta x^2)$. In contrast to \eqref{eq:4_4_4}, the error is not only due to the nonlinear part, since the discretisation error of the third derivative is in this case also in $\mathcal{O}(\Delta x^2)$. The results imply that the \texttt{mod1-central} and \texttt{mod2-central} schemes are AP for the $u$ component. If the time discretisation is GSA, then the schemes are AP (for all components). Numerically we observe that both schemes are third order accurate and AA which justifies the modifications. 

\subsection{Active Flux scheme with biased individual upwinding}

Although the \texttt{mod2-central} scheme delivers compelling results in the numerical experiments, we introduce an additional alternative scheme. It performs similarly to, if not better than, the \texttt{mod2-central} scheme, but shows that central discretisation is not the only viable approach. Inspired by the spatial discretisation of the HKdV equations presented in \cite{biswas_etal_2025}, we consider another choice of upwinding $F_2$ and define
\begin{align*}
	F_{2,\mathrm{biased}}^+ = \begin{pmatrix}
		0 & 0 & 1 \\
		0 & -\frac{1}{2\varepsilon^2} & 0 \\
		0 & 0 & 0
	\end{pmatrix} && \text{and} && F_{2,\mathrm{biased}}^- = \begin{pmatrix}
		0 & 0 & 0 \\
		0 & -\frac{1}{2\varepsilon^2} & 0 \\
		\frac{1}{\varepsilon^2} & 0 & 0	\end{pmatrix}.
\end{align*}
$u$ is upwinded from the left, $\psi$ from the right, and $\rho$ centrally.
Numerical experiments indicate that using the standard $D^\pm$ operators in the point value update equation of $\rho$ results in a scheme that is only second order. However, using the modified operators $D_\mathrm{mod}^\pm$ increases the order of accuracy to three. Due to the imposed directional bias, which is not based on characteristics, we call the resulting scheme the Active Flux scheme with biased individual upwinding or short the \texttt{biased} scheme.

By conducting a formal asymptotic analysis for this scheme (discretised in time with a type ARS IMEX Runge-Kutta scheme), we obtain \eqref{eq:3_3_3}, \eqref{eq:3_3_5}, and \eqref{eq:3_3_6} from the update equation of the cell averages and
\begin{subequations}
\label{eq:4_5_1}
	\begin{align}
		u_{i+\frac{1}{2}}^{[k](0)} &= 
	\begin{aligned}[t]
		&u^n_{i+\frac{1}{2}} + \Delta t \sum_{j=1}^{k}\tilde a_{k+1,j} \left(-\max\left((u_{i+\frac{1}{2}}^{(0)[j-1]},0\right)D^+u_{i+\frac{1}{2}}^{(0)[j-1]}\right.\\
		&\left.-\min\left(u_{i+\frac{1}{2}}^{(0)[j-1]},0\right)D^-u_{i+\frac{1}{2}}^{(0)[j-1]} \right) +\Delta t \sum_{j=1}^k a_{k,j}\left( - D^+\psi_{i+\frac{1}{2}}^{(0)[j]} \right) ,
	\end{aligned}\\
	u_{i+\frac{1}{2}}^{n+1(0)} &= \begin{aligned}[t]
		&u^n_{i+\frac{1}{2}} + \Delta t \sum_{j=1}^{\sigma}\tilde b_j \left(-\max\left(u_{i+\frac{1}{2}}^{(0)[j-1]},0\right)D^+u_{i+\frac{1}{2}}^{(0)[j-1]}\right.\\
		&\left.-\min\left(u_{i+\frac{1}{2}}^{(0)[j-1]},0\right)D^-u_{i+\frac{1}{2}}^{(0)[j-1]} \right) +\Delta t \sum_{j=1}^s b_{j}\left( - D^+\psi_{i+\frac{1}{2}}^{(0)[j]}  \right) 
	\end{aligned}
	\end{align}
\end{subequations}
and
\begin{equation}
	\psi_{i+\frac{1}{2}}^{(0)[j]} = \frac{1}{2} \left(D^+_\mathrm{mod}\rho_{i+\frac{1}{2}}^{(0)[j]} + D^-_\mathrm{mod}\rho_{i+\frac{1}{2}}^{(0)[j]} \right)
	\label{eq:4_5_2}
\end{equation}
and
\begin{equation}
	 \rho_{i+\frac{1}{2}}^{(0)[j]} = D^-u_{i+\frac{1}{2}}^{(0)[j]}
	 \label{eq:4_5_3}
\end{equation}
for all $j = 1, \dots, s$ from the update equations of the point values. Discretising in space with the \texttt{biased} scheme does not yield any unnatural conditions of the kind of \eqref{eq:unnat}. As before, it follows that \eqref{eq:3_3_3} and \eqref{eq:4_5_1} are temporal discretisations of 
\begin{equation}
	\begin{aligned}
	\frac{\mathrm{d}}{\mathrm{d}t}U_i^{(0)}(t) = \begin{aligned}[t]
			& -\frac{1}{\Delta x}\left(\frac{1}{2}\left(u_{i+\frac{1}{2}}^{(0)}(t)\right)^2 -\frac{1}{2}\left(u_{i-\frac{1}{2}}^{(0)}(t)\right)^2  \right)  \\
			&- \frac{1}{\Delta x^3}\left(\left(3U_i^{(0)[j]}  - 5u_{i+\frac{1}{2}}^{(0)[j]} + 4u_{i+\frac{3}{2}}^{(0)[j]} - 3U_{i+2}^{(0)[j]} + u_{i+\frac{5}{2}}^{(0)[j]}\right)\right.\\
			&\left.- \left(3U_{i-1}^{(0)[j]}  - 5u_{i-\frac{1}{2}}^{(0)[j]} + 4u_{i+\frac{1}{2}}^{(0)[j]} - 3U_{i+1}^{(0)[j]} + u_{i+\frac{3}{2}}^{(0)[j]}\right)\right).
		\end{aligned}
	\end{aligned}
	\label{eq:4_5_5}
\end{equation}
and
\begin{equation}
	\begin{aligned}
		\frac{\mathrm{d}}{\mathrm{d}t} u_{i+\frac{1}{2}}^{(0)}(t) = \begin{aligned}[t]
			& -\bigg(\max\left(u_{i+\frac{1}{2}}^{(0)}(t) ,0\right) \frac{1}{\Delta x}\left(2u_{i-\frac{1}{2}}^{(0)}(t) - 6U_i^{(0)}(t) + 4 u_{i+\frac{1}{2}}^{(0)}(t)\right) \\
			& + \min\left(u_{i+\frac{1}{2}}^{(0)}(t) ,0\right) \frac{1}{\Delta x}\left(-4u_{i+\frac{1}{2}}^{(0)}(t) + 6 U_{i+1}^{(0)}(t) - 2 u_{i+\frac{3}{2}}^{(0)}(t)  \right) \bigg)  \\
			&-\frac{1}{\Delta x^3} \left(6 U_{i-1}^{(0)[j]} -34u_{i-\frac{1}{2}}^{(0)[j]} + 48 U_i^{(0)[j]} -42 U_{i+1}^{(0)[j]}\right. \\
			&\left.+30 u_{i+\frac{3}{2}}^{(0)[j]}  -12U_{i+2}^{(0)[j]}+4u_{i+\frac{5}{2}}^{(0)[j]}\right).	
	\end{aligned}
	\end{aligned}
	\label{eq:4_5_4}
\end{equation}
In addition, if the IMEX Runge-Kutta scheme is GSA, we obtain from \eqref{eq:3_3_5}, \eqref{eq:3_3_6}, \eqref{eq:4_5_2}, and \eqref{eq:4_5_3},  
\begin{equation}
	\begin{aligned}
		\Psi_i^{(0)}(t) = \begin{aligned}[t]
		 &\frac{1}{\Delta x^2}\left(\left(-4u_{i+\frac{1}{2}}^{(0)}(t) + 6U_{i+1}^{(0)}(t) -2 u_{i+\frac{3}{2}}^{(0)}(t) \right)\right. \\
		 &\left.- \left(-4u_{i-\frac{1}{2}}^{(0)}(t) + 6U_{i}^{(0)}(t) -2 u_{i+\frac{1}{2}}^{(0)}(t) \right) \right)
		\end{aligned}
	\end{aligned}
	\label{eq:4_5_11} 
\end{equation}
and \eqref{eq:3_3_13} as well as 
\begin{equation}
	 \rho_{i+\frac{1}{2}}^{(0)}(t) = \frac{1}{\Delta x}\left(-4u_{i+\frac{1}{2}}^{(0)}(t) + 6U_{i+1}^{(0)}(t) -2 u_{i+\frac{3}{2}}^{(0)}(t) \right)
	 \label{eq:4_5_9}
\end{equation}
and
\begin{equation}
	\psi_{i+\frac{1}{2}}^{(0)}(t) = \frac{1}{\Delta x^2}\left(3U_i^{(0)}(t)  - 5u_{i+\frac{1}{2}}^{(0)}(t) + 4u_{i+\frac{3}{2}}^{(0)}(t) - 3U_{i+2}^{(0)}(t) + u_{i+\frac{5}{2}}^{(0)}(t)\right).
	\label{eq:4_5_10}
\end{equation}
Equation \eqref{eq:4_5_4} is a discretisation of \eqref{eq:th_1} with an error in $\mathcal{O}(\Delta x)$ and \eqref{eq:4_5_5} is a discretisation of \eqref{eq:th_2} with an error in $\mathcal{O}(\Delta x^2)$, implying that the \texttt{biased} scheme is AP for the component $u$. If the time discretisation scheme is GSA, then the scheme is AP since \eqref{eq:3_3_13}, \eqref{eq:4_5_11}, \eqref{eq:4_5_9}, and \eqref{eq:4_5_10} are discretisations of \eqref{eq:th_4}, \eqref{eq:th_6}, \eqref{eq:th_3}, and \eqref{eq:th_5} with an error in $\mathcal{O}(\Delta x^2)$. Our numerical experiments demonstrate that the \texttt{biased} scheme achieves third order accuracy and is AA. This shows that a much finer theoretical analysis of the order of accuracy is required.

\section{Numerical experiments}
\label{sec:num_exp}
The formal asymptotic analysis indicates that a GSA scheme is the preferable choice for temporal discretisation. Thus, unless otherwise specified, we use the $\mathrm{ARS}(4,4,3)$ scheme \cite{ascher_etal_1997}, which is GSA. The coefficients of this scheme are
\begin{align*}
	A = \begin{pmatrix}
		\frac{1}{2} & 0 & 0 & 0 \\[2pt]
		\frac{1}{6} & \frac{1}{2} & 0 & 0 \\[2pt]
		-\frac{1}{2} & \frac{1}{2} & \frac{1}{2} & 0 \\[2pt]
		\frac{3}{2} & - \frac{3}{2} & \frac{1}{2} & \frac{1}{2}
	\end{pmatrix} && \text{and} && b = \begin{pmatrix}
		\frac{3}{2} \\[2pt]
		 - \frac{3}{2} \\[2pt]
		  \frac{1}{2} \\[2pt]
		   \frac{1}{2}
	\end{pmatrix}
\end{align*}
and
\begin{align*}
	\tilde A = \begin{pmatrix}
		0 & 0 & 0 & 0 & 0\\[2pt]
		\frac{1}{2} & 0 & 0 & 0 & 0\\[2pt]
		\frac{11}{18} & \frac{1}{18} & 0 & 0 & 0 \\[2pt]
		\frac{5}{6} & -\frac{5}{6} & \frac{1}{2} & 0 & 0 \\[2pt]
		\frac{1}{4} & \frac{7}{4} & \frac{3}{4} & -\frac{7}{4} & 0
	\end{pmatrix} && \text{and} && b = \begin{pmatrix}
		\frac{1}{4} \\[2pt]
		\frac{7}{4} \\[2pt]
		\frac{3}{4} \\[2pt]
		-\frac{7}{4} \\[2pt]
		0
	\end{pmatrix}.
\end{align*}
The schemes are implemented in Julia \cite{julia_2017} and the results are visualised with Makie.jl \cite{makie}.

\begin{figure}[t]
	\centering
	\includegraphics[scale=0.6]{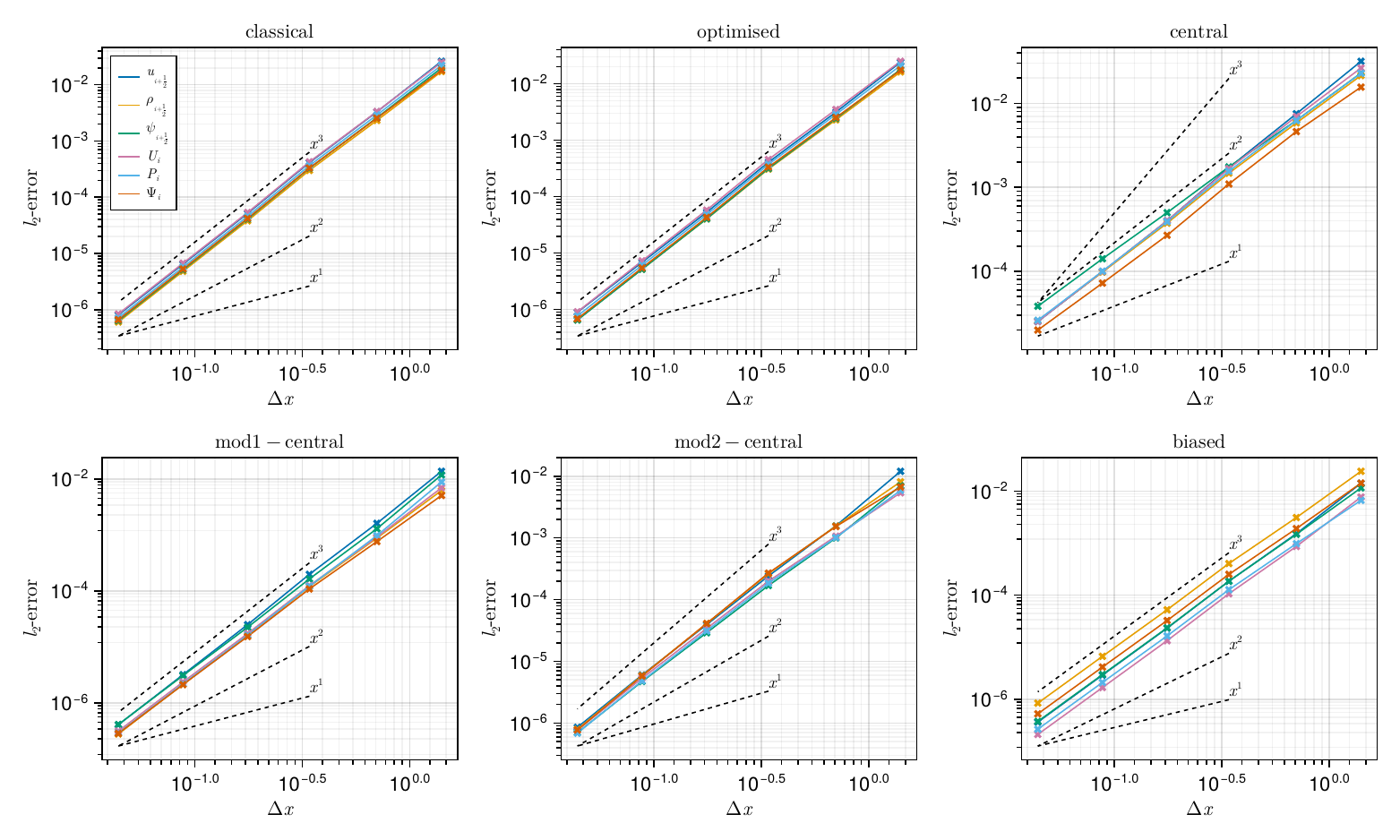}
	\caption{Experiment 1: The $l_2$-error of the point values and cell averages measured with respect to the spatial step size $\Delta x$ with $\varepsilon = 0.8$ for the different IMEX Active Flux schemes. All schemes are third-order accurate except the \texttt{central} scheme which is only second-order accurate.}
	\label{fig:eoc1}
\end{figure}

\subsection{Experiment 1: Travelling wave solution of the HKdV equations}
This experiment investigates the ability of our schemes to solve an initial value problem for the HKdV equations, i.e.\ a value of $\varepsilon$ close to $1$. For this experiment, $\varepsilon = 0.8$ is chosen. As initial data, we use the travelling wave solution of the HKdV equation that tends to the single-soliton solution of the KdV equation and its derivatives with $c=\frac{1}{2}$ and $x_0=0$ at $t=0$. Since it is not analytically available, the travelling wave solution is computed numerically as described in \cite{biswas_etal_2025} on the interval $[x_{\mathrm{min}}, x_{\mathrm{max}}] = [-47, 47]$ with periodic boundary conditions.
In contrast to \cite{biswas_etal_2025}, eighth-order central finite difference formulas \cite{fornberg_1988} are used instead of Fourier spectral differentiation. We find that using high-order finite difference formulas and a sufficiently fine grid (42000 points) is sufficient for the experiments in this text. The numerical solution at the grid points is interpolated with cubic splines to obtain a continuous function $\tilde q_\mathrm{ref}: [x_{\mathrm{min}}, x_{\mathrm{max}}]  \rightarrow \R^3$. Shifting this function to the right with respect to time, i.e.\ $q_\mathrm{ref}(t,x) = \tilde q_\mathrm{ref}(x-ct)$, gives an accurate approximation to the travelling wave solution.

Given the initial data and reference solution, we compute the numerical solutions on the interval $[x_{\mathrm{min}}, x_{\mathrm{max}}] = [-45, 45]$ with periodic boundary conditions using the Active Flux schemes proposed in this work. For all Active Flux schemes we choose $C_\mathrm{CFL} = 0.35$ if not stated otherwise.

We assess the experimental convergence rate, the error on a given grid and for long times, stability, and energy conservation.
First, we compute the numerical solution at $t=4$ on grids with $64$, $128$, $256$, $512$, $1024$, and $2048$ cells and measure the $l_2$-error for all point values and cell averages, see Figure \ref{fig:eoc1}. The \texttt{central} scheme achieves an order of accuracy of about two, while the other schemes show third order accuracy.

In order to make a comparison between the schemes on the same coarse grid with $64$ cells, the maximum $l_2$-error is computed as
	\begin{align*}
		\mathcal{E}_{l_2,\mathrm{max}}(\Delta x) =\max\left( \max_{v \in \Gamma_1} \mathcal{E}_{l_2}(v,\Delta x), \max_{V \in \Gamma_2} \mathcal{E}_{l_2}(V,\Delta x) \right),
	\end{align*}
where in case of the HKdV equations $\Gamma_1 = \{u,\rho, \psi\}$ and $\Gamma_2 = \{U,P, \Psi\}$. The maximum $l_2$-error of the different schemes is displayed in Table \ref{tab:comp_e1}. The \texttt{mod2-central} scheme admits the lowest error.

\begin{table}[t]
\renewcommand{\arraystretch}{1.5}
\centering
\begin{tabular}{|M{2.2cm}|M{2.2cm}|M{2.2cm}|M{2.2cm}|M{2.2cm}|M{2.2cm}|N}
\hline
      \texttt{classical}    &     \texttt{optimised} &   \texttt{central} &  \texttt{mod1-central} &  \texttt{mod2-central} &  \texttt{biased} &    \\
 \hline
      \multicolumn{7}{c}{\text{maximum $l_2$-error on a grid with $64$ cells at $t=4$}} \\
 \hline
   0.0262      &   0.0251       &0.0317&  0.0137 & \textbf{0.0120} &0.0240 & \\
 \hline
\multicolumn{7}{c}{\text{maximum $l_2$-error on a grid with $512$ cells at $t=180$}} \\
 \hline
   0.0116  &   0.0121     &  0.0094 & 0.0019 & 0.0015 &\textbf{0.0006} &\\
 \hline
 \multicolumn{7}{c}{\text{maximum $C_\mathrm{CFL} \leq 1$ on a grid with $512$ cells at $t=4$}} \\
 \hline
   \textbf{1.000}    &    0.450     &  0.450 & 0.450 & 0.450 &\textbf{1.000} &\\
 \hline
  \multicolumn{7}{c}{\text{conservation defect with $ \Delta x = \frac{90}{256}, \frac{90}{512}, \frac{90}{1048}$ (top to bottom) at $t=360$ }} \\
 \hline
 $1.0000\cdot10^{-2}$ &    $1.0423\cdot10^{-2}$   & $7.6617\cdot 10^{-3}$ &  $1.6431\cdot10^{-3}$ & $1.3863 \cdot 10^{-3}$& $\mathbf{4.1771 \cdot{10^{-4}}}$ &\\
 \hline
    $1.2609 \cdot 10^{-3}$    &  $1.3176\cdot10^{-3}$  & $8.6421 \cdot 10^{-4}$ & $2.0559 \cdot 10^{-4}$ & $1.6120 \cdot 10^{-4}$ & $\mathbf{5.5125\cdot10^{-5}}$ &\\
 \hline
   $1.5779 \cdot 10^{-4}$     & $1.6476 \cdot10^{-4}$   & $9.5892\cdot10^{-5}$ & $2.5729\cdot10^{-5}$ & $1.9574 \cdot  10^{-5}$ &$\mathbf{6.9489\cdot10^{-6}}$&\\
 \hline
\end{tabular}
\caption{Results of Experiment 1 performed for the different Active Flux schemes discretised in time with the $\mathrm{ARS}(4,4,3)$ scheme. The best performance is highlighted in boldface. }
\label{tab:comp_e1}
\end{table}

To investigate the long-time accuracy of the methods, we compute the maximum $l_2$-error on a grid with $512$ cells at the time at which the wave has travelled once over the computational domain $[-45, 45]$, i.e.\ at $t = 180$ for $c = \frac{1}{2}$ (see Table \ref{tab:comp_e1}). The \texttt{biased} scheme has the lowest error for this long simulation time.

In order to experimentally investigate stability bounds of the different schemes, we determine the maximum $l_2$-error of the numerical solution at $t=4$ on a grid with $512$ for different values of $C_\mathrm{CFL} \leq 1$. We take extreme error growth as a sign of instability. Table \ref{tab:comp_e1} shows the maximum value of $C_\mathrm{CFL}$ for which the error is not higher than ten times the error for $C_\mathrm{CFL} = 0.1$. The \texttt{biased} scheme and the \texttt{classical} scheme are the only schemes for which this holds for all tested values of $C_\mathrm{CFL} \leq 1$.

Lastly, we investigate the conservation of energy by the different numerical schemes. It is not clear if a given Active Flux scheme conserves a discrete energy and how this energy is defined. We compare the proposed schemes by measuring
\begin{align*}
	E_\mathrm{H,dis}(\Delta x, t) = \Delta x \frac{1}{2} \sum_{i=1}^N \left((u_{i+\frac{1}{2}}(t))^2 + \varepsilon^2 (\rho_{i+\frac{1}{2}}(t))^2 + \varepsilon^2 (\psi_{i+\frac{1}{2}}(t))^2 \right),
\end{align*}
which can be considered a discrete counterpart of the energy of the HKdV equations. $N$ denotes the number of cells. We define the conservation defect at time $t$ as 
\begin{align*}
	\Delta E_\mathrm{H,dis}(\Delta x, t) =  1 - \frac{E_\mathrm{H,dis}(\Delta x, t)}{E_\mathrm{H,dis}(\Delta x, 0)}.
\end{align*}
This value is displayed in Table \ref{tab:comp_e1}, where we chose $t = 360$ and $\Delta x = \frac{90}{256}, \frac{90}{512}, \frac{90}{1048}$. All schemes conserve the energy satisfactorily but there is a factor of 20 in the defect between the best and the worst scheme on the same grid. The \texttt{biased} scheme delivers the best results followed by the \texttt{mod2-central} scheme. 

\begin{figure}[t]
	\centering
	\includegraphics[scale=0.6]{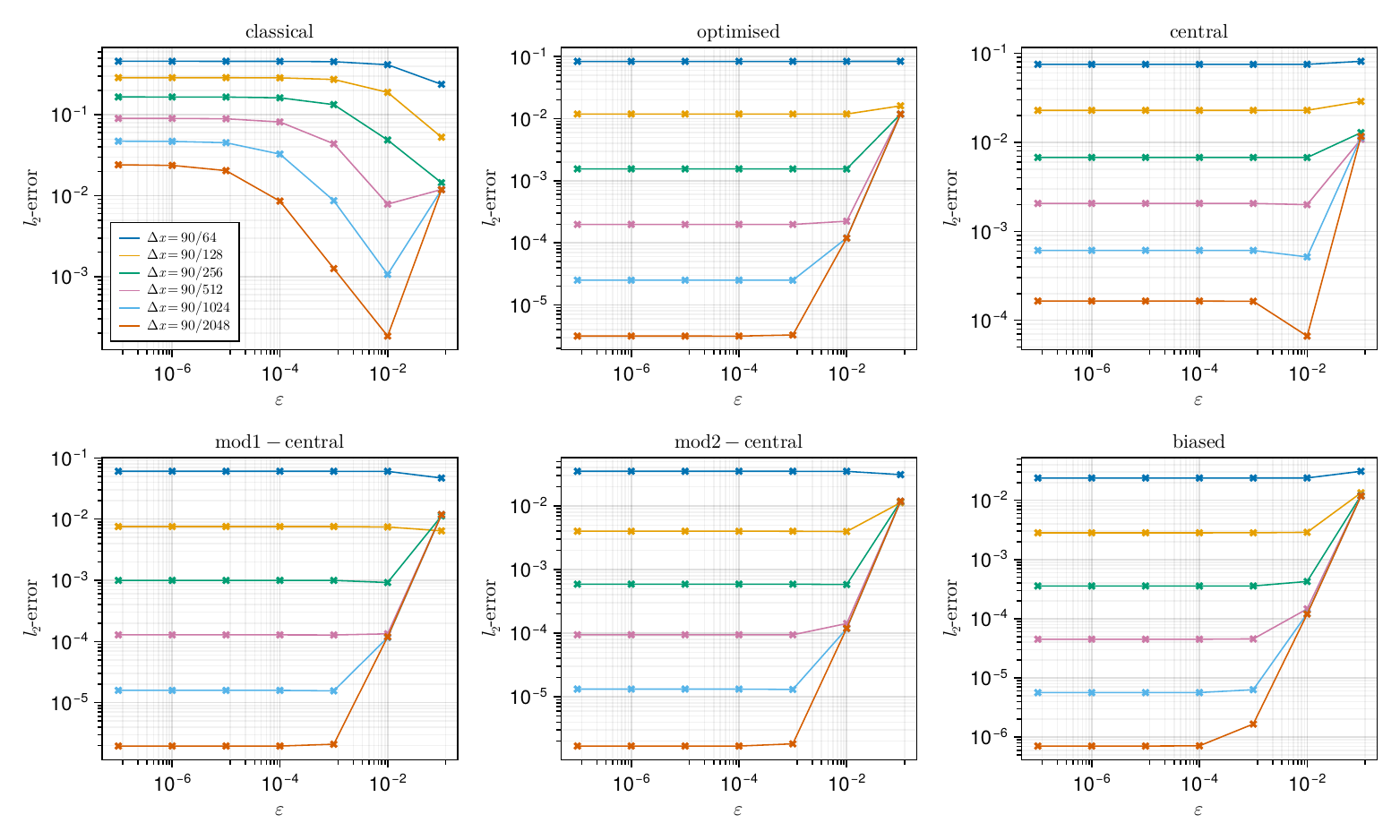}
	\caption{Experiment 2: The $l_2$-error of $u_{i+\frac{1}{2}}$ with respect to the relaxation parameter $\varepsilon$ for the different IMEX Active Flux schemes and on different grids.}
	\label{fig:ap}
\end{figure}

Observe that the results obtained for the \texttt{classical} and \texttt{optimised} schemes are similar, which is not unexpected since the schemes are closely related. The results of the experiment conducted for the \texttt{combined} scheme discretised in time with the explicit SSPRK3 scheme are also comparable to the results obtained with the  \texttt{classical} and \texttt{optimised} scheme: it achieves third order accuracy, the maximum $l_2$-error on a coarse grid is $0.0246$ and for long times $0.0118$, the maximum value of $C_\mathrm{CFL}$ is $0.5$, and the conservation defects are $1.0123 \cdot 10^{-2}$ ($256$ cells), $1.2834 \cdot 10^{-3}$ ($512$ cells), and $1.6085 \cdot 10^{-4}$ ($1024$ cells). 

\subsection{Experiment 2: Single-soliton solution of the KdV equation}

\begin{figure}[p]
	\centering
	\includegraphics[scale=0.6]{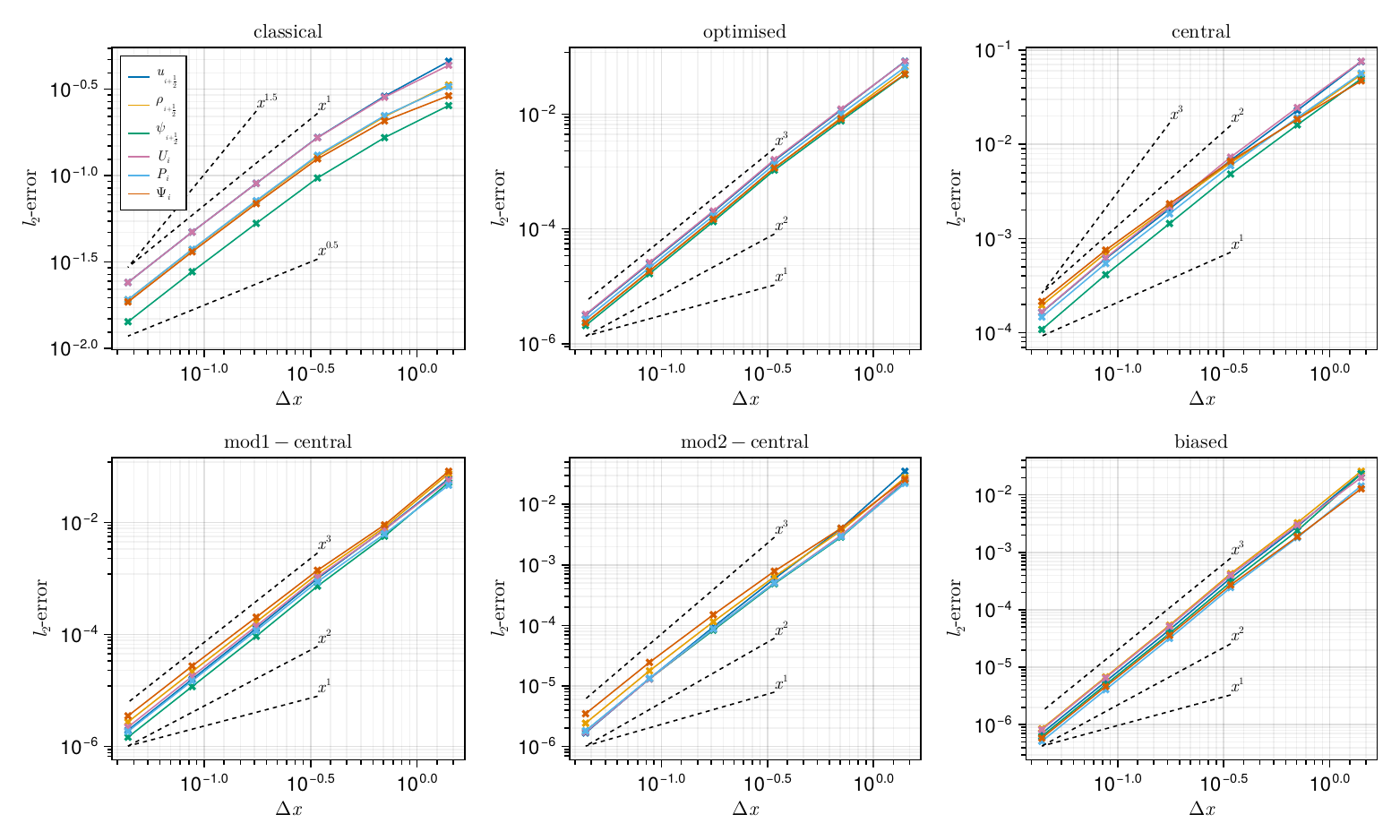}
	\caption{Experiment 2: The $l_2$-error of the point values and cell averages with respect to the spatial step size $\Delta x$ for $\varepsilon = 10^{-7}$. The \texttt{classical} scheme only achieves first-order accuracy, the \texttt{central} scheme second-order accuracy, and the other schemes third-order accuracy.}
	\label{fig:eoc2}
\end{figure}
\begin{figure}[p]
	\centering
	\includegraphics[scale=0.6]{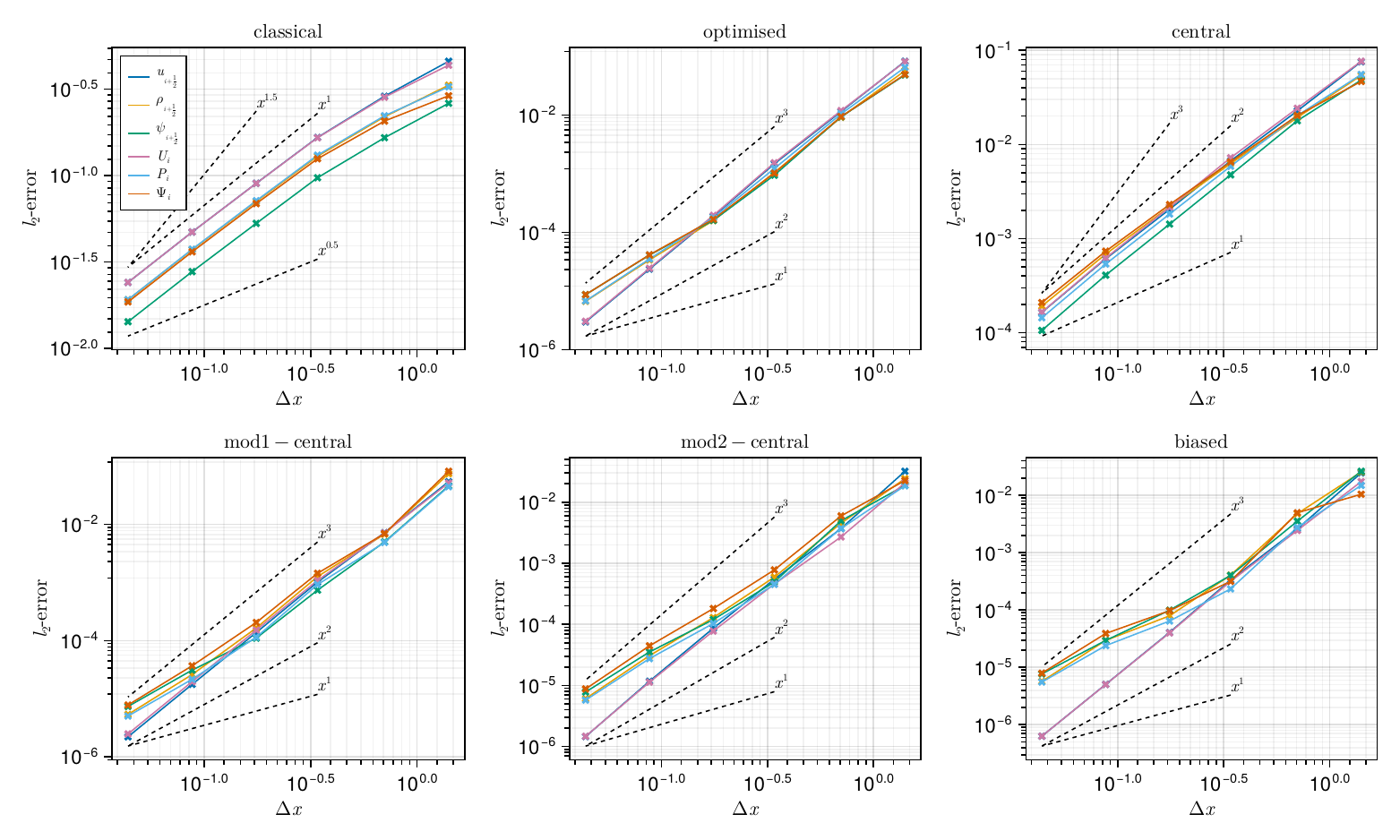}
	\caption{Experiment 2: The $l_2$-error of the point values and cell averages with respect to the spatial step size $\Delta x$ for the different Active Flux schemes discretised in time with the $\mathrm{ARS}(2,3,3)$ scheme.}
	\label{fig:ars233}
\end{figure}

Next, we investigate the ability of the schemes to compute numerical approximations to the limit equations of the HKdV equations, i.e.\ we study the behaviour of the methods in the limit $\epsilon \to 0$. The interval $[x_{\mathrm{min}}, x_{\mathrm{max}}] = [-45, 45]$ is chosen, periodic boundary conditions are imposed, and the CFL number $C_\mathrm{CFL}=0.35$ is used. The single-soliton solution \eqref{eq:soliton} of the KdV equations (with $x_0=0$ and $c= \frac{1}{2}$) and its first two derivatives serve as reference solution. In particular, the initial data are
\begin{align*}
	u_0(x) &= 3c \, \mathrm{sech}^2\left(\frac{\sqrt{c}}{2}x\right), \\
	\rho_0(x) &= -3 c^\frac{3}{2}\mathrm{sech}^2\left(\frac{\sqrt{c}}{2}x\right) \mathrm{tanh}\left(\frac{\sqrt{c}}{2}x\right), \\
	\psi_0(x) &= 3 c^2 \tanh ^2\left(\frac{\sqrt{c}}{2}x\right) \text{sech}^2\left(\frac{\sqrt{c}}{2}x\right)-\frac{3}{2} c^2 \text{sech}^4\left(\frac{\sqrt{c}}{2}x\right).
\end{align*}
The errors are computed with respect to this reference solution and not with respect to the solution of the HKdV equations for the given initial value problem. This decision is justified since we want to determine if the scheme can solve an initial value problem for the limit equations.

In Figure \ref{fig:ap} the $l_2$-error of the numerical solution at $t=4$ is shown for different values of $\varepsilon$ and different spatial step sizes. For simplicity only the $l_2$-error of the point values $u_{i+\frac{1}{2}}$ is displayed. The errors of the other point values and cell averages behave similarly. For all schemes, the error becomes constant for a fixed $\Delta x$ as $\varepsilon \rightarrow 0$. This implies that for the smallest value of $\varepsilon$ considered, i.e.\ $\varepsilon = 10^{-7}$, the measured error is the discretisation error of the limit scheme (cf. Figure \ref{fig:ap_prop}). Observe that if $\Delta x \rightarrow 0$, the error decreases. This is shown in more detail in Figure \ref{fig:eoc2}, where the error of all point values and cell averages is displayed for fixed $\varepsilon = 10^{-7}$ and different values of $\Delta x$. Figures \ref{fig:ap} and \ref{fig:eoc2} are numerical evidence that the schemes are AP, since, as $\varepsilon \rightarrow 0$, the schemes become a (consistent and stable) discretisation of the limit scheme. For a scheme that is not AP, the measured error for small $\epsilon$ would remain constant upon grid refinement, for example.

The accuracy of the \texttt{classical} scheme drops to first-order, while all the other schemes retain the same order of accuracy as for large $\varepsilon$. Consequently, the numerical experiment suggests that all schemes except the \texttt{classical} scheme are AA.

\begin{table}[t]
\renewcommand{\arraystretch}{1.5}
\centering
\begin{tabular}{|M{2.2cm}|M{2.2cm}|M{2.2cm}|M{2.2cm}|M{2.2cm}|M{2.2cm}|N}
\hline
      \texttt{classical}    &     \texttt{optimised} &   \texttt{central} &  \texttt{mod1-central} &  \texttt{mod2-central} &  \texttt{biased} &    \\
 \hline
      \multicolumn{7}{c}{\text{maximum $l_2$-error on a grid with $64$ cells at $t=4$}} \\
 \hline
   0.4576      &   0.0837   &0.0760&  0.0820 & 0.0350 &\textbf{0.0258} & \\
 \hline
\multicolumn{7}{c}{\text{maximum $l_2$-error on a grid with $512$ cells at $t=180$}} \\
 \hline
   3.5441   &   0.0578     &  0.3014 & 0.0190 & \textbf{0.0086} &0.0123 &\\
 \hline
 \multicolumn{7}{c}{\text{maximum $C_\mathrm{CFL} \leq 1$ on a grid with $512$ cells at $t=4$}} \\
 \hline
   \textbf{1.000}    &    0.400     &  0.450 & 0.425 & 0.425 &\textbf{1.000} &\\
 \hline
  \multicolumn{7}{c}{\text{conservation defect with $ \Delta x = \frac{90}{256}, \frac{90}{512}, \frac{90}{1048}$ (top to bottom) at $t=360$ }} \\
 \hline
 $7.1943\cdot10^{-1}$ &    $3.2949\cdot10^{-2}$   & $1.0313\cdot 10^{-1}$ &  $-9.5614\cdot10^{-3}$ & $\mathbf{5.0983\cdot 10^{-3}}$& $7.0442 \cdot{10^{-3}}$ &\\
 \hline
    $6.1156\cdot10^{-1}$    &  $4.3219\cdot10^{-3}$  & $2.4094\cdot 10^{-2}$ & $-1.3954 \cdot 10^{-3}$ & $\mathbf{6.5200 \cdot 10^{-4}}$ & $9.1454\cdot10^{-4}$ &\\
 \hline
   $4.7749\cdot 10^{-1}$     & $5.4399 \cdot10^{-4}$   & $4.1704\cdot10^{-3}$ & $-1.8348\cdot10^{-4}$ & $\mathbf{8.7434 \cdot  10^{-5}}$ &$1.1553 \cdot10^{-4}$&\\
 \hline
\end{tabular}
\caption{Results of Experiment 2 performed for the different IMEX Active Flux scheme discretised in time with the $\mathrm{ARS}(4,4,3)$ scheme. The best performance is highlighted in boldface.}
\label{tab:comp_e2}

\end{table}

\begin{remark}[Importance of the GSA property]
Figure \ref{fig:ars233} shows the $l_2$-error of the point values and cell averages for different spatial step sizes if the schemes are discretised in time with the $\mathrm{ARS}(2,3,3)$ scheme, which is not GSA and has the following coefficients: 
\begin{align*}
	A = \begin{pmatrix}
		\gamma & 0 \\[1pt]
		1-2\gamma & \gamma
	\end{pmatrix} && \text{and} && b = \begin{pmatrix}
		\frac{1}{2} \\[1pt]
		\frac{1}{2}
	\end{pmatrix}
\end{align*}
and
\begin{align*}
	\tilde A = \begin{pmatrix}
		0 & 0 & 0 \\[1pt]
		\gamma & 0 & 0 \\[1pt]
		\gamma -1 & 2(1-\gamma) & 0
	\end{pmatrix} && \text{and} && \tilde b = \begin{pmatrix}
		0 \\[1pt]
		\frac{1}{2} \\[1pt]
		\frac{1}{2}
	\end{pmatrix}
\end{align*}
with $\gamma = \frac{3+\sqrt{3}}{6}$. The results obtained for the \texttt{classical} and \texttt{central} scheme are similar to those shown in Figure \ref{fig:eoc2}: the \texttt{classical} scheme achieves first order accuracy and the \texttt{central} scheme second order accuracy. The other schemes display third-order accuracy for $u$, but only a reduced order for $\rho$ and $\psi$, implying that the schemes are AP but with exception of the central scheme are not AA for all components. 

In addition, the matrices $A$ of the $\mathrm{ARS}(2,3,3)$ and $\mathrm{ARS}(4,4,3)$ schemes have the same entries on the diagonal. This allows us to precompute the LU-decomposition of the resulting matrices in the time discretisation. Hence, even though the $\mathrm{ARS}(4,4,3)$ scheme has more stages than the $\mathrm{ARS}(2,3,3)$ scheme, their efficiency is almost the same since computing the LU-decomposition is the most expensive part. This justifies our choice to employ the $\mathrm{ARS}(4,4,3)$ scheme to discretise in time.
\end{remark}

For the remainder of the experiment $\varepsilon = 10^{-7}$ is chosen. As in the first experiment, to compare the different schemes, we measure their maximum $l_2$-error of the numerical solutions at $t=4$ computed on a coarse gird ($64$ cells) and at $t=180$ on a grid with $512$ cells. The results are shown in Table \ref{tab:comp_e2}. The \texttt{biased} scheme has the smallest error on the coarse grid and the \texttt{mod2-central} scheme the lowest for long times. Additionally, we determine the highest values of $C_\mathrm{CFL}$, for which the maximum error of a numerical solution computed on a grid with $512$ cells at $t=4$ is smaller than ten times the error with $C_\mathrm{CFL} = 0.1$. As shown in Table \ref{tab:comp_e2}, only the \texttt{classical} and \texttt{biased} scheme achieve this for all $C_\mathrm{CFL} \leq 1$.

Lastly, we investigate the conservation of a discrete counterpart of the energy conserved by the KdV equation. We define 
\begin{align*}
			E_\mathrm{dis}(\Delta x, t) = \Delta x \sum_{i=1}^N \frac{1}{2}\left(u_{i+\frac{1}{2}}(t)\right)^2
\end{align*}
and the conservation defect 
\begin{align*}
	\Delta E_\mathrm{H}(\Delta x, t) =  1 - \frac{E_\mathrm{H}(\Delta x, t)}{E_\mathrm{H,dis}(\Delta x, 0)}.
\end{align*}
The defect at $t=360$ for $\Delta x = \frac{90}{256}, \frac{90}{512}, \frac{90}{1048}$ is shown in Table \ref{tab:comp_e2}. There are 2 orders of magnitude in the measured defect between the best and the worst scheme. The \texttt{mod2-central} scheme gives the best results, followed by the \texttt{biased} scheme. Observe that the discrete energy of the numerical solution computed with the \texttt{mod1-central} scheme increases and that the defect of the classical nearly stagnates upon mesh refinement.
 \begin{figure}[t]
 	\includegraphics[scale=0.6]{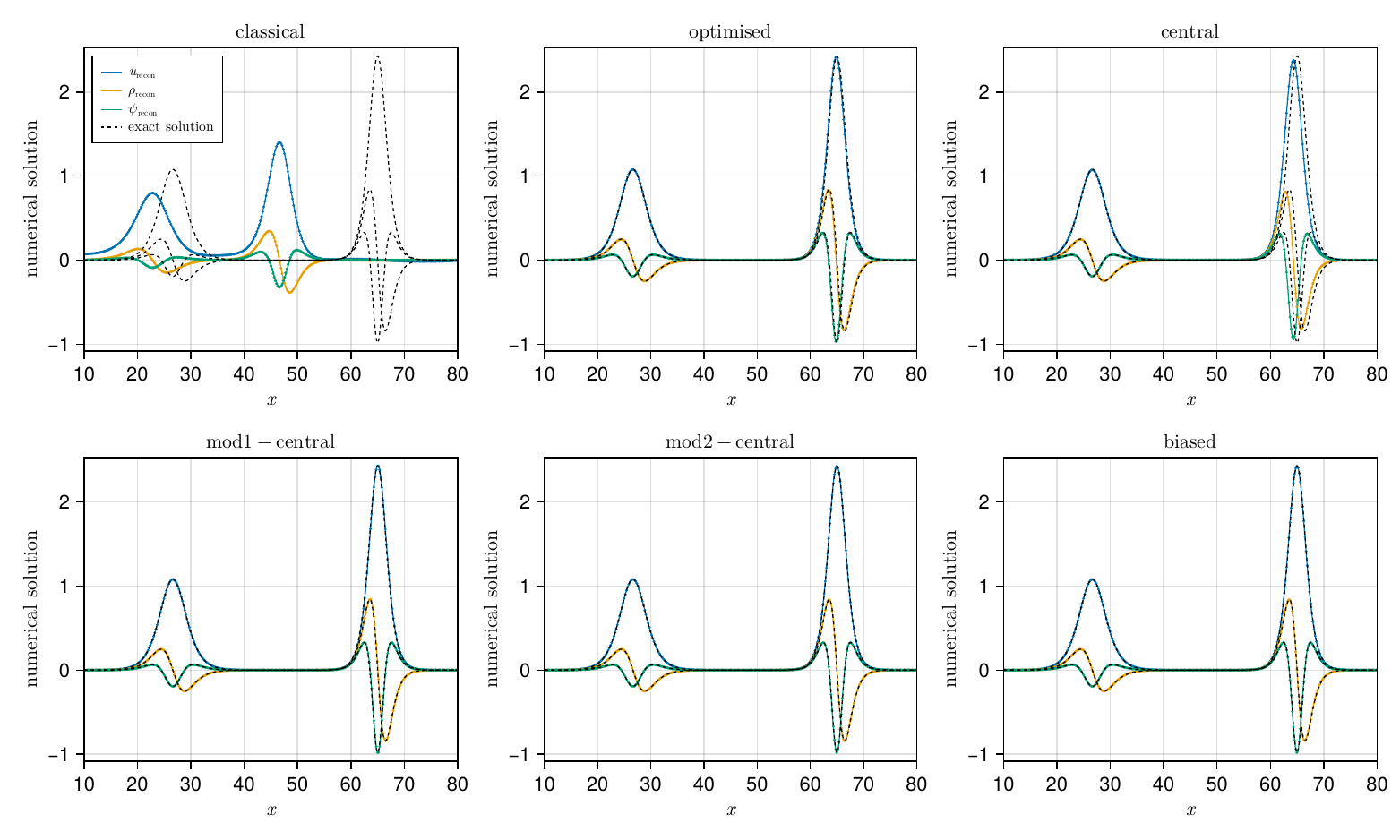}
	\caption{Experiment 3: The numerical double-solution solution computed with the different Active Flux schemes in comparison to the exact solution at $t = 100$. The faster soliton has already overtaken the slower one.}
	\label{fig:doub_sol}
\end{figure}

 \subsection{Experiment 3: Double-soliton solution of the KdV equation}

\begin{figure}[t]
	\centering
	\includegraphics[scale=0.6]{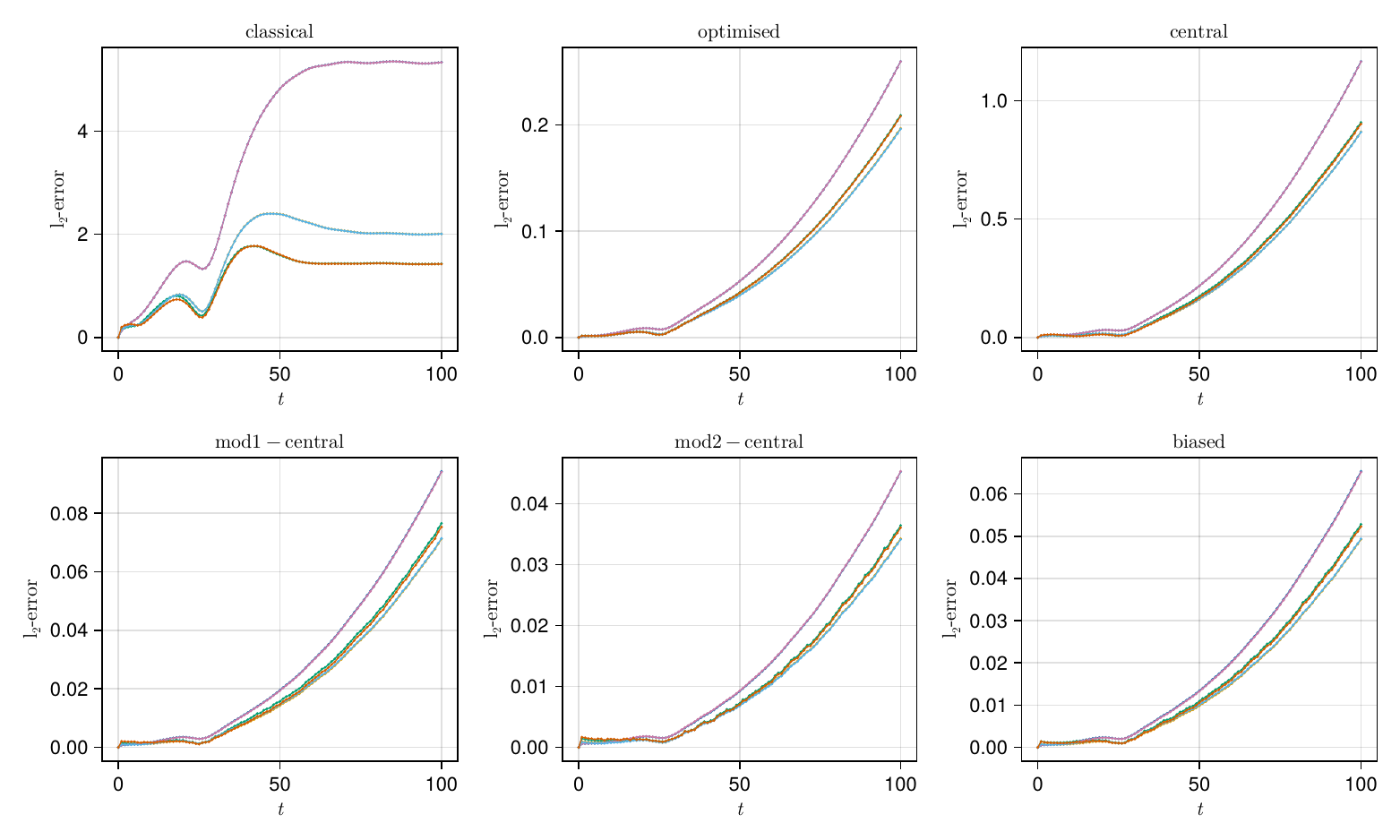}
	\caption{Experiment 3: The $l_2$-error of the numerical double-solution solution computed with the different Active Flux schemes as a function of time $t$.}
	\label{fig:doub_sol_2}
\end{figure}
 In this experiment, we are employing the different Active Flux schemes to numerically compute the double-soliton solution of the KdV equation. We set $\varepsilon = 10^{-7}$ and $C_\mathrm{CFL} = 0.35$ for all parts of the experiment, choose $[x_{\mathrm{min}}, x_{\mathrm{max}}] = [-30.0, 80.0]$ divided into $512$ cells, and impose periodic boundary conditions. As initial data,
\begin{align*}
	u_0(x) &= 2(\delta_2^2 - \delta_1^2) \frac{\delta_2^2 \mathrm{cosech}(\gamma_2(0,x))^2 +\delta_1^2 \mathrm{sech}(\gamma_1(0,x))^2}{(\delta_2\mathrm{coth}(\gamma_2(0,x)) +\delta_1\mathrm{tanh}(\gamma_1(0,x)))^2}, \\
	\rho_0(x) &= Du_0(x), \\
	\psi_0(x) &= D\rho_0(x)
\end{align*}
is used, where $\gamma_1$ and $\gamma_2$ are defined as in Theorem \ref{th:doub_sol} and $D$ is a fourth order central finite difference operator \cite{fornberg_1988} with spatial step size $h = 0.0005$; $\delta_1 = 0.3$ and $\delta_2 = 0.45$ are selected. The initial data consist of two solitons, a smaller, slower one that is located on the right of a faster, larger one. The analytical solution of this initial value problem for the limit equations is the function given in Theorem \ref{th:doub_sol} and its first and second derivative. This solution acts as the reference solution in this experiment. The first and second derivatives of it are again computed using fourth order central finite differences.

Figure \ref{fig:doub_sol} shows the numerical solution computed with the different Active Flux schemes at $t = 100$. At this time, the faster soliton has already overtaken the slower soliton. The $l_2$-error of the numerical solution for various times $t$ is displayed in Figure \ref{fig:doub_sol_2}. The \texttt{mod2-central} scheme is the most accurate out of all the different schemes, while the \texttt{classical} and the \texttt{central} schemes yield solutions that deviate strongly from the exact one.

\section{Conclusion and outlook}
\label{sec:con}
After presenting the \texttt{combined} scheme, i.e.\ the standard semi-discrete Active Flux scheme for the HKdV equations, we derived a closely related scheme, the \texttt{classical} scheme,  by individually upwinding the nonlinear and stiff part of the equations. This scheme allows the discretisation in time with IMEX Runge-Kutta schemes. Formal asymptotic analysis was conducted for this scheme using an IMEX Runge-Kutta scheme of type ARS. The theoretical findings show that the fully-discrete scheme is AP for the $u$ component of the HKdV equations. If the time discretisation is globally stiffly accurate (GSA), then the scheme is AP for all components. This is confirmed experimentally, yet the scheme is only first-order accurate for small $\varepsilon$.

To find schemes with superior properties, we considered different ways to upwind the stiff part of the Jacobian. First, the \texttt{optimised} scheme was obtained by modifying and rescaling particular elements of the upwinding matrices of the stiff operator. In addition, using central discretisation for the terms discretised in time implicitly, the \texttt{central} scheme was developed. Since the scheme turned out to achieve only second-order accuracy, two modifications were introduced, the \texttt{mod1-central} and \texttt{mod2-central} schemes. Both of them achieve an order of accuracy of three. Lastly, to demonstrate an alternative to central discretisation, we presented the \texttt{biased} scheme, where we imposed predefined upwind directions in the upwinding matrices of the stiff operator. Formal asymptotic analysis for these schemes gave the same results as for the \texttt{classical} scheme, i.e.\ that they all are AP.

The numerical experiments conducted in this work verify that all of the schemes presented here are AP, which is generally not trivial to achieve. In fact, the experiments show that they already appear to be AP if they are discretised in time with a scheme that is not GSA. All schemes apart from the \texttt{classical} scheme are AA, a property clear from the experiments but incompletely understood theoretically. We observed this property (with the exception of the \texttt{central} scheme) only if we used a GSA time discretisation. By comparing the accuracy, maximum CFL number, and conservation of a discrete energy, the \texttt{mod2-central} scheme and the \texttt{biased} scheme emerge as the schemes with the best properties.

Further analysis of the presented alternative point value update equations and resulting schemes is a compelling topic for future research into the reasons for their superior properties compared to the \texttt{classical} scheme. One option is to implement the limit schemes and measure the rate at which the numerical approximations of the Active Flux schemes tend to the numerical approximations of the corresponding limit schemes. Moreover, it would be interesting to further study the role of the time-discretisation by using different IMEX Runge-Kutta schemes of type I and type II and comparing the results. It is also worth identifying simple criteria for a given Active Flux scheme to be AP, without having to redo the asymptotic analysis. In the future, we will apply the IMEX Active Flux schemes developed here to more complex models of dispersive free-surface waves. The applicability of IMEX schemes is not restricted to hyperbolised high-order PDEs, and we will explore their use in other contexts.

\subsection*{Acknowledgements}
Christian Klingenberg acknowledges funding by the Deutsche Forschungsgemeinschaft (DFG, German Research Foundation) within SPP 2410 Hyperbolic Balance Laws in Fluid Mechanics: Complexity, Scales, Randomness (CoScaRa), project No.\ 525941602.

\end{document}